\documentclass[11pt]{article}

\author{Bing-Long Chen\\[8pt]}
\title{\textbf{Local  regularity  of spacetimes under Ricci curvature and Lie derivative conditions\\
}}
\date{Oct. 4, 2025}
\usepackage{latexsym}
\usepackage{amsmath}
\usepackage{amssymb,amscd}
\usepackage[mathscr]{eucal}
\newtheorem{thm}{Theorem}[section]

\newtheorem{cor}[thm]{Corollary}

\newtheorem{prop}[thm]{Proposition}

\newtheorem{defn}{Definition}[section]
\newtheorem{rem}{Remark}[section]
\newtheorem{claim}{Claim}
\numberwithin{equation}{section}
\newenvironment{pf}{{\noindent \it  Proof.}}{{\hfill$\Box$}\\}

\usepackage{amsmath}
\usepackage{amssymb}
\usepackage{bm}
\usepackage[T1]{fontenc} 

\begin{document}

\maketitle
\let\thefootnote\relax\footnotetext{AMS Mathematics Subject Classification Numbers: Primary
53c50; Secondary 83c20. }

\begin{abstract}
	In this paper,  we derive a general  regularity estimate  for  \textbf{any}  4-d spacetime,  in terms of  a priori  bounds of the Ricci curvature  and  Lie derivative of the Lorentzian metric with respect  to an \textbf{arbitrary}  timelike vector field. 	\end{abstract}

  \section{Introduction} 
  
  According to the theory of general relativity,  the  gravity can be   described by a $4$-dimensional Lorentzian manifold $(\bar{M}, \bar{g})$,  where  $\bar{g}$ is a Lorentzian metric,   satisfying  the Einstein's equation:
  \begin{equation} \label{e1.1}
  \bar{R}ic-{1}/ {2} \bar{R} \bar{g}=\kappa T
  \end{equation}
  where  $\kappa$ is a  constant,  $T$ is the  energy-momentum tensor determined  by the distribution of  matter  or  energy   in the universe.  
  
  After specifying the equations governing the physical fields which make up the  energy-momentum tensor $T$,  and imposing a harmonic  gauge condition,  (\ref{e1.1}) can be reduced to a  second order hyperbolic system.  Given an initial data set obeying constraint equations on a spacelike hypersurface,   the Cauchy problem to   (\ref{e1.1})  can be solved (\cite{CB}).    The dynamical behavior of the  solutions to the Cauchy problem,  is a central  issue   in the   mathematical theory of  general relativity.   
   
     On the other hand,  from geometric perspective, we consider the following problem.   Suppose we have a piece of spacetime  $(\bar{M}, \bar{g})$,  without prescribing a Cauchy  surface, what can we say about the regularity of $(\bar{M}, \bar{g})$ around a point $x_0\in \bar{M}$.  Here, the regularity refers to the existence of a ``coordinate system'' covering a large portion of the region around $x_0$ on which the coefficients  of metric $\bar{g}$ have  bounded norms in certain functional space.  It is instructive to note  that from the example of plane waves (\S 5.9 in \cite{HE}),  even in the vacuum case, if no further conditions are imposed, it is impossible to derive any regularity estimate.      
  
   There are many reasons to introduce an auxiliary timelike vector field to the picture (whose  existence  is guaranteed if $\bar{M}$ is time-orientable). For instance, solving the Cauchy problem will produce time functions whose gradient vector fields are timelike.   
    
   The purpose of the paper  is to    investigate the regularity of a  4-dimensional  spacetime $(\bar{M},\bar{g})$ in the presence of  an \textbf{arbitrarily} given timelike vector field $X$ on $\bar{M}$.  
               
         To present the main result,  we introduce   a Riemannian metric $\hat{g}=\hat{g}_{X}$ canonically associated to $X$:  
   \begin{equation} \label{1.2}
   \hat{g}=2 u^{-2} X_{\ast}\otimes X_{\ast}+\bar{g}
   \end{equation}
   where $u^{2}=-\bar{g}(X,X)>0$ and $X_{\ast}$ is the co-vector field dual to $X$,  i.e. $X_{\ast}$ satisfies $X_{\ast}(Y)=\bar{g}(X,Y)$ for any vector field $Y$.  By  using $\hat{g}$, one can measure the lengths of all tensors on $\bar{M}$.

         One of the main results of the paper is the following:  
     \begin{thm} \label{t1.1}   
    Let  $(\bar{M}, \bar{g})$ be  a four dimensional spacetime with a timelike vector field $X$.  Suppose  the closure of  the $\hat{g}$-geodesic ball ${B}_{\hat{g}}(x_0, a)$ of radius $a>0$ is compact and   \begin{equation}\begin{split} \label{1.3}
&\sup_{{B}_{\hat{g}}(x_0, a)} |\bar{R}ic|_{\hat{g}}  \leq  a^{-2}\\
& \sup_{{B}_{\hat{g}}(x_0, a)} u^{-1}|\mathcal{L}_{X}\bar{g}|_{\hat{g}} \leq a^{-1}, 
  \end{split}
\end{equation}
where $\bar{R}ic$ is the Ricci curvature of $\bar{g}$, $\mathcal{L}_{X}\bar{g}$ is the  Lie derivative,   $|\cdot|_{\hat{g}}$ is the norm defined by the  Riemannian metric  $\hat{g}$. 

 Then for any $\epsilon>0$, there exist constants $c_1>0, c_2>0$ depending only on $\epsilon$ so that for any $p> 0$, there is   a smooth immersion   $\Psi$  from $\{z\in \mathbb{R}^4: |z|< c^{-1}_1a\}$ to $\bar{M}$ with   $\Psi(0)=x_0$ such that:\\
 i)   the  image of $\Psi$ covers  $B_{\hat{g}}(x_0,  c_2^{-1} a)$; \\
 ii) the coordinates $\{z^{\alpha}\}$ are $\hat{g}=\Psi^{\ast} \hat{g}$-harmonic, i.e. $\triangle_{\hat{g}} z^{\alpha}=0$;\\
  iii)  $\bar{g}=\Psi^{\ast} \bar{g}$  satisfies the following estimates: 
 \begin{equation}
 \begin{split}  \label{1.4}
& \sup_{\{|z|<c^{-1}_1a\}}|\bar{g}_{\alpha\beta}-  \eta_{\alpha\beta}|(z)    \leq \epsilon\\
  &  \frac{1}{a^4}\int_{\{|z|<c^{-1}_1a\}} |\partial \bar{g}_{\alpha \beta}|^p  dz   \leq \epsilon a^{-p},
 \end{split}
  \end{equation}
where     $(\eta_{\alpha \beta})=\text{diag} (-1,+1,+1,+1)$  is the Minkowski metric,  $\bar{g}_{\alpha\beta}=\bar{g}(\partial_{z^\alpha},\partial_{z^{\beta}})$.     \end{thm}

 In Theorem \ref{t1.1},  we  assume that the closure of  the $\hat{g}$-geodesic ball ${B}_{\hat{g}}(x_0, a)$ is compact.  Otherwise, there already  were  some pathologies  for  the spacetime in the ball  ${B}_{\hat{g}}(x_0, a)$, which should  be  avoided in the  discussions of this paper.  

Roughly speaking,   Theorem \ref{t1.1} says that the $W^{1,p}$ estimate  can be obtained provided that  the Ricci curvature and  Lie derivative of the Lorentzian metric are bounded.

A  by-product of proving Theorem \ref{t1.1} is  the $W^{1,p}$ estimate for the timelike  vector field  $X$.  It should be noted that  the Lie derivative  condition   (\ref{1.3})  is only involving  the symmetric part of  derivative $D_{\bar{g}}X$ of $X$. 

 \begin{cor} \label{c1.2} Under  the conditions  of Theorem \ref{t1.1},   let   $X=X^{\alpha}\frac{\partial}{\partial z^{\alpha}}$ in  the coordinates $\{z^{\gamma}\}$,      we have 
\begin{equation}
 \begin{split}  \label{1.5}
& \sup_{\{|z|<c^{-1}_1a\}}|u^{-1}X^{\alpha}-\delta_{\alpha0}|(z)  \leq \epsilon\\
  & \frac{1}{a^4}\int_{\{|z|<c^{-1}_1a\}} u^{-p}|\partial_{\alpha} X^{\beta}|^{p}dz    \leq \epsilon a^{-p},  \end{split}
  \end{equation}
where $\delta_{\alpha 0}$ is equal to 1 for $\alpha=0$,  and 0 otherwise. 
 \end{cor}

  By Sobolev imbedding theorem,  for any  $0<\alpha<1$, the $C^{\alpha}$ norms (in coordinates $\{z^{\alpha}\}$ of Theorem \ref{t1.1}) of  $\bar{g}$, $X$ and  $\hat{g}$ are uniformly bounded.    
  
In Theorem \ref{t1.1}, if we assume additionally the covariant derivative of $\mathcal{L}_{X}\bar{g}$ is bounded, i.e.
\begin{equation}\begin{split} \label{1.6}
\sup_{\hat{B}(x_0, a)}  u^{-1} |D_{\bar{g}}(\mathcal{L}_{X}\bar{g})|_{\hat{g}}  \leq a^{-2},   
\end{split}
\end{equation}
 one can obtain $W^{2,p}$ estimate (for any $p>0$): 
\begin{equation} \label{1.7}
 \frac{1}{a^4}\int_{\{|z|<c^{-1}_1a\}} |\partial^2 \bar{g}_{\alpha \beta}|^p+u^{-p}|\partial^2 X^{\alpha}|^pdz  \leq c_p a^{-2p},
 \end{equation}
where  $c_p$ is a  positive constant depending only on $\epsilon$ and $p$. 
 
  In particular,  (\ref{1.7})  gives an intrinsic and coordinate-free curvature estimate:
 \begin{equation} 
  		\begin{split}\label{1.8}
  			\frac{1}{vol_{\hat{g}}(\hat{B}(x_0,\frac{a}{2}))}\int_{\hat{B}(x_0,\frac{a}{2})}|\bar{R}m(\bar{g})|_{\hat{g}}^p dv_{\hat{g}} \leq c_p a^{-2p}. 
			\end{split}
  	\end{equation}     
	
	Higher order regularity of $\bar{g}$ and $X$ can also be obtained, if  the correspondingly more  derivatives  of  the Ricci curvature and Lie derivative are bounded. In  particular,   the  $L^{\infty}$ bound of the  curvature tensor of $\bar{g}$ can be controlled:      
	
	    \begin{cor} \label{c1.3}  If we replace the condition (\ref{1.3}) in Theorem \ref{t1.1} with  the following 
  \begin{equation}\begin{split} \label{1.9}
& \sup_{\hat{B}(x_0, a)} |\bar{R}ic|_{\hat{g}}+a|D_{\bar{g}}\bar{R}ic |_{\hat{g}}  \leq  a^{-2}\\
& \sup_{\hat{B}(x_0, a)} u^{-1}|\mathcal{L}_{X}\bar{g}|_{\hat{g}}+  a u^{-1} |D_{\bar{g}}(\mathcal{L}_{X}\bar{g})|_{\hat{g}}+a^2  u^{-1} |D^2_{\bar{g}}(\mathcal{L}_{X}\bar{g})|_{\hat{g}} \leq a^{-1},
\end{split}
\end{equation}
we have  \begin{equation}
 \begin{split} \label{1.10}
 & \sup_{B_{\hat{g}}(x_0,\frac{a}{2})}|\bar{R}m|_{\hat{g}} \leq C a^{-2},\\
\end{split}
  \end{equation}
  where $C$ is a universal constant. 
  \end{cor}

 Corollary \ref{c1.3} is an improvement of Theorem 0.1 in \cite{A1}  and Theorem 1.3 in  \cite{Chen1}.   
  If $\bar{R}ic\equiv 0$, and $X$ is a timelike Killing field, i.e. $\mathcal{L}_{X}\bar{g}=0$,    the conditions in (\ref{1.9}) hold automatically.   If in addition, $(\bar{M}, \bar{g})$ is  geodesically complete, one can show that $(\bar{M}, \hat{g})$ is also geodesically complete (Theorem 3.3 in \cite{Chen1}). By taking  $a\rightarrow \infty$ in (\ref{1.9}) and (\ref{1.10}), we have   $\bar{R}m \equiv 0$, i.e. $(\bar{M}, \bar{g})$ is  flat. See \cite{A1}    \cite{Lich55} \cite{EP} for these  pioneering  works.

It is worthy to be   mentioned  that when the  dimension of the spacetime  $\geq 5$, Theorem 1.1  does  not  hold.  A counterexample can be constructed as follows.  Take $(N, g_{N})$ to be a complete Ricci flat and non-flat Riemannian manifold of dimension $\geq 4$, then  $(\bar{M},\bar{g})=(\mathbb{R}\times N, -dt^2+g_{N})$ is a complete Ricci flat and non-flat spacetime of dimension $\geq 5$ with an obvious timelike Killing field $X={\partial}/ {\partial t}$.


  \begin{prop}  \label{p1.4}
  	Under the conditions  of  Theorem \ref{t1.1}, if we further assume \begin{equation}
  		\begin{split}\label{1.11}
  		vol_{\hat{g}}(\hat{B}(x_0,\frac{a}{2})) \geq v_0 a^4 
  		\end{split}
  		\end{equation}
  	for some $v_0>0$,  and if we allow the  constants $c_1, c_2, c_p$ to depend on  $v_0$, then the map $\Psi$ in Theorem \ref{t1.1} can be chosen to be a diffeomorphism.   
	 \end{prop}
	
	In other words, under (\ref{1.11}),   $\{z^{\alpha}\}$ can  be a  coordinate system in Theorem \ref{t1.1} in the ordinary sense.  Combining Theorem \ref{t1.1} and Proposition \ref{p1.4},  one can establish a compactness theorem which is useful when we need to take limits from a sequence  of  Lorentzian manifolds.

We point out  some  connections or comparisons   of  Theorem \ref{t1.1} with  some   breakdown criterions  for vacuum Einstein equations  and  bounded $L^2$ curvature conjecture.

   1)   The coordinate system in Theorem \ref{t1.1} is an immersed coordinate system, that means we allow the volume is collapsed.    
    
  2)       In many works  on vacuum Einstein equations,  in doing energy estimate,  the vanishing of  covariant derivatives of the Ricci curvature are usually substantially employed.   In Theorem \ref{t1.1},   the Ricci curvature is only assumed to be bounded.     
  
   3)  As we emphasized before,  Theorem \ref{t1.1} can handle the situation without assuming the existence of a Cauchy surface,  in contrast to  the classical works on vacuum Einstein equations making  uses of the  foliation of maximal or constant mean curvature  surfaces, cf. \cite{A2}\cite{KR10}\cite{W} \cite{KRS}.

   Theorem \ref{t1.1} and Corollary \ref{c1.2} can be applied to relativistic (ideal) fluids, where the regularity of the velocity field of the fluid is also mainly concerned.   Recall that the relativistic (ideal) fluids satisfy  \begin{equation}
	\begin{split}
		\label{1.12}
		& \bar{R}_{\alpha\beta}-\frac{1}{2} \bar{R} \bar{g}_{\alpha\beta}=\kappa [(\mu+p) u_{\alpha}u_{\beta}+p \bar{g}_{\alpha\beta}]
	\end{split}
\end{equation}
where  $u^{\alpha}$ is the  unit 4-velocity vector field of  particles in the fluid,   $\mu$, $p$  are the energy and pressure functions respectively.   
Apparently,  $u^{\alpha}$ can be   a natural reference frame. The $\hat{g}$-bound of  $\bar{R}ic$ is equivalent to the bound of $\mu$ and $p$. One can also choose   $X^{\alpha}= f u^{\alpha}$, where $f$ is  some  thermodynamical  function of the fluid,  provided that a state equation describing the details of the fluid  is  given (Chapter 2 in \cite{Lich67}).    Theorem \ref{t1.1} and Corollary \ref{c1.2} can also be applied to many other   circumstances where timelike vector fields can be naturally singled out.

 Theorem \ref{t1.1} can be generalized  by replacing   the $L^{\infty}$ norms  in condition (\ref{1.3}) with normalized  $L^{p}$ norms for  scales $\leq a$ ($p>2$ for Ricci curvature  and $p>4$ for Lie derivative),  see  (\ref{7.7}) (\ref{7.8}) and Theorem \ref{t7.1}.

  We discuss the  strategies of  proving Theorem \ref{t1.1} and the arrangement of the paper.  In Section 2, we  formulate  the result in Theorem \ref{t1.1} to an  estimate of  certain weak harmonic radius.  In Section 3, we calculate the equations of $\bar{g}, \hat{g}, X$ which are needed for later arguments in Sections 4 and 6.   The overall strategy is an argument by contradiction. Suppose Theorem \ref{t1.1}  does  not  hold,   there is  a sequence of marked and normalized  Lorentzian manifolds $(\bar{M}_j, \bar{g}_j, \hat{g}_j, X_j, \bar{x}_j)$ violating Theorem \ref{t1.1} such that the  weak harmonic radius $r^{w,h}(W^{1,p_0}, \epsilon_0)$ at marked point $x_j$ is 1 and  greater than 1/2 in a large  $\hat{g}_j-$ball of radius $r_j \rightarrow \infty$, and \begin{equation} \nonumber
\sup_{{B}_{\hat{g}_j}(\bar{x}_j, r_j)} |{R}ic(\bar{g}_j)|_{\hat{g}_j}  +u_j^{-1}|\mathcal{L}_{X_j}\bar{g}_j|_{\hat{g}_j} \rightarrow 0.
\end{equation} 

The naive idea is the following. If one could extract a smooth limit (from a  subsequence) of  geodesically complete vacuum spacetime with a timelike Killing field,  the limit must be flat by \cite{Chen1}.  The point is that  the flat limit has infinite large weak harmonic radius, this will contradict with the normalization condition that the harmonic radius is 1 at the origin.  There are several difficulties in this argument. First of all, the sequence $({B}_{\hat{g}_j}(\bar{x}_j, r_j), \hat{g}_j, \bar{x}_j)$ is possibly volume-collapsed. The Gromov-Hausdorff limit, if it exists, is possibly a singular and lower dimensional space.  This prevents the direct application of  \cite{Chen1} to the limit.  So one need an  appropriate  convergent limit to maintain some smooth structure for which we have the vacuum Einstein equation and timelike Killing field.  In each coordinate chart $\Omega_i$ given by weak harmonic radius condition $r^{w,h}(W^{1,p_0}, \epsilon_0) \geq 1/2$,  a subsequence of  $(\bar{g}_j, \hat{g}_j, X_j)$ will converge to a limit $(\bar{g}_{\infty}, \hat{g}_{\infty}, X_{\infty})$  which is  vacuum with  a timelike Killing field.  To apply the global argument of \cite{Chen1},    the  point is to consider the convergence of the transition functions  between different charts.  This motivates us to consider  groupoids.   We only use the idea of  groupoids, instead of developing in detail the theory of Lorentizian groupoids.  In Section 4,   we  prove that after choosing subsequences,  in each coordinate chart $\Omega_i$,  the limit $(\bar{g}_{\infty}, \hat{g}_{\infty}, X_{\infty})$ is flat.  This enables us to smooth  the Riemannian metrics $\hat{g}_j$  a little so that the resulting   metrics $(\hat{g}_j)_{\epsilon}$  have  very small curvatures (Section 5), i.e., $|Rm((\hat{g}_j)_{\epsilon})|\rightarrow 0$.  Pulling back  $\bar{g}_j, X_j$  to large balls of  tangent spaces (using the exponential maps of $(\hat{g}_j)_{\epsilon}$),   the closeness of $\bar{g}_j, X_j$  with the flat limit  will ultimately imply a very large weak harmonic radius at the origin (Section 6), which gives a contradiction.   Proposition \ref{p1.4} is proved in Section 7.

 \section{Gromov-Hausdorff convergence}
  We shall formulate the result in Theorem \ref{t1.1} in terms of certain  weak harmonic radius, the new feature is that the timelike vector field $X$ has to be incorporated in the Lorentzian setting. 
  \begin{defn} \label{d2.1}
  For any $\epsilon>0$, $p>0$, $P\in \bar{M}$,  we define the weak $W^{1,p}$-harmonic radius of accuracy $\epsilon$ at $P$ (denoted by  $r^{w,h}(W^{1,p}, \epsilon)(P)$) to be the supremum of all $r>0$ such that  there is a smooth  immersion  $\Psi$ defined on a neighborhood of $\{z\in \mathbb{R}^4: |z|\leq r\}$ with $\Psi(0)=P$ such that $z^{\alpha}$ is harmonic with respect to the pulled back Riemannian metric $\hat{g}\triangleq \Psi^{\ast} \hat{g}$, i.e. 
   \begin{equation} \label{2.1}
  	\begin{split}
  		\triangle_{\hat{g}}z^{\alpha}=0, \ \ \ \  \alpha=0,1,2,3,
  	\end{split}
  \end{equation}
and 
  \begin{equation} \label{2.2}
  	\begin{split}
  		\sup_{\{|z|\leq r\}} |\hat{g}_{\alpha\beta}-\delta_{\alpha\beta}|+ |\bar{g}_{\alpha\beta}-\eta_{\alpha\beta}|+|u^{-1}X^{\alpha}-\delta_{\alpha 0}|\leq \epsilon
  		\end{split}
  	\end{equation}
  
   \begin{equation} \label{2.3}
  		 	r^{p-4}\int_{\{|z|\leq r\}} |\partial \hat{g}_{\alpha\beta}|^p+ |\partial \bar{g}_{\alpha\beta}|^p+u^{-p}|\partial X^{\alpha}|^pdz \leq \epsilon.	
  \end{equation}
  
  \end{defn}
  
  \begin{rem} \label{r2.1}

 Similarly, one can define the weak $C^{\alpha}$-harmonic radius $r^{w,h}(C^{\alpha},\epsilon)$ of accuracy $\epsilon$ as well.  Here, we call the quantity   ``weak'' harmonic radius because the map $\Psi$ is only an immersion.  If $\Psi$ is required to be a diffeomorphism, we can just call it  a $W^{1,p}$ or $C^{\alpha}$-harmonic radius as usual, and denoted it by $r^{h}(W^{1,p}, \epsilon)(P)$ or $r^{h}(C^{\alpha},\epsilon)(P)$,  depending on which norms are chosen. 
\end{rem}

 Theorem \ref{t1.1} and Corollary \ref{c1.2} can be rephrased simply as 
 
 \begin{equation} \label{2.4}
 		r^{w,h}(W^{1,p}, \epsilon)(x_0)\geq C^{-1} a
 		\end{equation}
 	 for some constant $C$ depending only on $\epsilon$ and  $p$. 
  
  Apparently,  it suffices  to prove the following equivalent and scaling invariant statement: 
 
 \begin{equation} \label{2.5}
 	\sup_{x\in  B_{\hat{g}}(x_0, a)}	d_{\hat{g}}(x, \partial B_{\hat{g}}(x_0, a)) 	[r^{w,h}(W^{1,p}, \epsilon)(x)]^{-1}\leq C.	
 	\end{equation}

By H$\ddot{o}$lder inequality, we only need to prove (\ref{2.5}) for large $p$. 

 To prove (\ref{2.5}),   we shall try to use a contradiction argument.

 By scaling the metric, we may  assume $a=1$.  Suppose  (\ref{2.5}) is not true, there are  $\epsilon_0>0$, $p_0>4$ and a sequence of pointed Lorentzian manifolds with timelike vector fields $(\bar{M}_j, \bar{g}_j, X_j, P_j)$  such that  

 \begin{equation}\begin{split} \label{2.6}
\sup_{{B}_{\hat{g}_j}(P_j, 1)} |{R}ic(\bar{g}_j)|_{\hat{g}_j}  \leq  1, \ \ \ \ 
\sup_{{B}_{\hat{g}_j}(P_j, 1)} u_j^{-1}|\mathcal{L}_{X_j}\bar{g}_j|_{\hat{g}_j} \leq 1, 
  \end{split}
\end{equation}

but  \begin{equation} \label{2.7}
 	\sup_{x\in  B_{\hat{g}_j}(P_j, 1)}	d_{\hat{g}_j}(x, \partial B_{\hat{g}_j}(P_j, 1)) 	[r^{w,h}(W^{1,p_0}, \epsilon_0)(x)]^{-1} \triangleq m_j\geq j\rightarrow \infty.	
 	\end{equation}
 Let $\bar{x}_j$ be a point in  $B_{\hat{g}_j}(P_j, 1)$ such that 
\begin{equation} \label{2.8}
 	d_{\hat{g}_j}(\bar{x}_j, \partial B_{\hat{g}_j}(P_j, 1)) 	[r^{w,h}(W^{1,p_0}, \epsilon_0)(\bar{x}_j)]^{-1}=m_j.	
 	\end{equation}
	Denote $\epsilon_j= r^{w,h}(W^{1,p_0}, \epsilon_0)(\bar{x}_j)$. Then $\epsilon^{-1}_j \geq m_j \rightarrow \infty$. Scale  $\bar{g}_j$ with factor $\epsilon_j^{-2}$, scale  $X_j$ so that $|X_j|_{\hat{g}_j}^{2}(\bar{x}_j)=1$, i.e.  $u_j(\bar{x}_j)=1$,   still denote them  with the  same notations. 
	For the new metric, we have  $B_{\hat{g}_j}(\bar{x}_j, m_j)\subset B_{\hat{g}_j}({P}_j, \epsilon^{-1}_j)$,
	\begin{equation} \label{2.9}\begin{split} r^{w,h}(W^{1,p_0}, \epsilon_0)(\bar{x}_j)& =1,\\ 
	r^{w,h}(W^{1,p_0}, \epsilon_0)({x})& \geq \frac{1}{2}, \ \ \ \  \text{on} \ \   B_{\hat{g}_j}(\bar{x}_j, \frac{1}{2}m_j), 	\end{split}\end{equation}	
	and 	
 \begin{equation}\begin{split} \label{2.10}
\sup_{{B}_{\hat{g}_j}(\bar{x}_j, \frac{1}{2}m_j)} |{R}ic(\bar{g}_j)|_{\hat{g}_j}  \leq  \epsilon^2_j, \ \ \ \ 
\sup_{{B}_{\hat{g}_j}(\bar{x}_j, \frac{1}{2}m_j)} u_j^{-1}|\mathcal{L}_{X_j}\bar{g}_j|_{\hat{g}_j} \leq \epsilon_j.
  \end{split}
\end{equation}		

In the subsequent sections (Sections 3-6), we shall try to derive a contradiction by taking suitable “limits” from the subsequence of  $({B}_{\hat{g}_j}(\bar{x}_j, \frac{1}{2}m_j), \bar{g}_j, \hat{g}_j, X_j, \bar{x}_j)$ as $j\rightarrow \infty$.  Currently,  at least,  we can prove that the   pointed Riemannnian manifolds   $(({B}_{\hat{g}_j}(\bar{x}_j, \frac{1}{2}m_j), \hat{g}_j, \bar{x}_j)$  have a  Gromov-Hausdorff convergent subsequence.   

The point is to prove the following Claim: 

\begin{claim}:  \label{claim1}  For any $r>0$, $\epsilon>0$, there exist  constants  $N(r,\epsilon)>0$ and  $j_0=j_0(r,\epsilon)>0$, such that for any $j\geq j_0$, the ball ${B}_{\hat{g}_j}(P, r)\subset {B}_{\hat{g}_j}(\bar{x}_j, \frac{1}{2}m_j)$ can be covered by at most $N(r,\epsilon)$  balls of radius $\epsilon$. 
\end{claim}
\begin{pf}
  For any ball $B_{\hat{g}_j}(P, r) \subset {B}_{\hat{g}_j}(\bar{x}_j, \frac{1}{2}m_j) $.  Let  $N(P, r, \epsilon)$ be the maximum  number of points in $B_{\hat{g}_j}(P, r)$ with mutual distances $\geq \epsilon$. It is clear $B_{\hat{g}_j}(P, r)$ can be covered by at most $N(P, r, \epsilon)$ geodesic balls of radius  $\epsilon$.   Since   $r^{w,h}(W^{1,p_0},\epsilon_0)(P) \geq 1/2$, there is a map $\Psi: \{|z|\leq 1/2\}\rightarrow \bar{M}_j$ with $\Psi(0)=P$ such that (\ref{2.2}) holds.  First, we  consider the case for  $r\leq  1/3$.  Let $\{y_k\}_{k=1}^{k_0}\subset B_{\hat{g}_j}(P, r)$ be $k_0$-points with mutual distances $\geq \epsilon$, and $\{\tilde{y}_k\}_{k=1}^{l}\subset  \{|z|<5/12\}$ with $\Psi(\tilde{y}_k)=y_k$. It is clear $|\tilde{y}_k-\tilde{y}_{k^{\prime}}|\geq \frac{1}{2} \epsilon$ for $k\neq k^{\prime}$, which implies  $k_0\leq C \epsilon^{-n}$, $n=dim \bar{M}_j$. Hence   $N(P, r, \epsilon) \leq C\epsilon^{-n}$ for $r\leq 1/3$. One can take $N(r,\epsilon)=C\epsilon^{-n}$ when $r\leq 1/3$,  for some universal constant $C$.  
  
  Now we handle  the general  case for $r\geq 1/3$. By an induction argument, we will prove that there is a universal constant $C>0$ such that  $N(\frac{l}{6}, \frac{1}{6})\leq C^{l}$, for $l=1,2,\cdots$.  Assume that  we have proved $N(\frac{l}{6}, \frac{1}{6})\leq C^{l}$ for $l\leq l_0$.   Suppose  
  $B_{\hat{g}_j}(P, \frac{l_0+1}{6}) \subset {B}_{\hat{g}_j}(\bar{x}_j, \frac{1}{2}m_j)$. For any point $x\in B_{\hat{g}_j}(P, \frac{l_0+1}{6})$, there is a point $x^{\prime}\in B_{\hat{g}_j}(P, \frac{l_0}{6})$ with $d_{\hat{g}_j}(x,x^{\prime})\leq 1/6$. This implies  that $B_{\hat{g}_j}(P, \frac{l_0+1}{6})$ can be covered by at most $N(\frac{l_0}{6}, \frac{1}{6})$ balls of radius   $1/3$. By the previous step, each of these $1/3$-balls can be covered by $N(\frac{1}{3},\frac{1}{6}) \leq C$  balls of radius $1/6$. This gives $N(\frac{l_0+1}{6},\frac{1}{6}) \leq N(\frac{l_0}{6},\frac{1}{6}) N(\frac{1}{3},\frac{1}{6}) \leq C^{l_0+1}$. 
  
    Since $N(\frac{l}{6}, \epsilon) \leq N(\frac{l}{6},\frac{1}{6}) N(\frac{1}{6},\epsilon) \leq C C^{l} \epsilon^{-n}$, one can take 
    $N(r,\epsilon)=Ce^{Cr} \epsilon^{-n}$.  The proof  is completed.  
    \end{pf}
    
     By compactness theorem of Gromov (Chapter 10 in \cite{Pe}),  there exists  a pointed complete metric space $(L,d_{L},\bar{x})$ such that a subsequence of  $({B}_{\hat{g}_j}(\bar{x}_j, \frac{1}{2}m_j), \hat{g}_j, \bar{x}_j) $  Gromov-Hausdorff  converges to  $(L,d_{L},\bar{x})$.   For simplicity of  notations,  we may assume the subsequence is just the original sequence.

   \section{Einstein equations}
   In this section,   we shall  compute the equations satisfied by the metric $\hat{g}$ and the vector field $X$.  These equations are presented in a coordinate free manner.  Proposition \ref{p3.1} computes the Ricci curvature of $\hat{g}$,  which can be regarded as the equations  for  $\hat{g}$.   Propositions \ref{p3.2} and \ref{p3.3}  compute  the equations satisfied by $u$ and a  horizontal 2-form $\Lambda$, which can be regarded as the main equations for $X$.   These computations are straightforward and tedious,  they can be bypassed for the first reading.

    In the followings,  if some  indices $\alpha,\beta,$ $\cdots$ of  tensors appear,  they are understood to come from  an arbitrarily fixed  coordinates $\{z^{\alpha}\}$.

  First of all, we decompose the covariant derivative $D^{\bar{g}}X$ of $X$ into symmetric  and anti-symmetric parts.  
  
  For vector fields $Y$ and $Z$, we define  
 \begin{equation} \label{3.1}
 \langle D^{\bar{g}}_YX, Z \rangle\triangleq \bar{A}(Y,Z)=\frac{1}{2} (\mathcal{L}_{X}\bar{g})(Y,Z)+\frac{1}{2} A(Y,Z) 
 \end{equation}
 where $A(Y,Z)=\langle D^{\bar{g}}_YX, Z \rangle-\langle D^{\bar{g}}_ZX, Y \rangle$ is anti-symmetric in $Y$ and $Z$. 
   Let  \begin{equation} \label{3.2} P= I+u^{-2}\langle X, \cdot \rangle_{\bar{g}} X=I-u^{-2}\langle X, \cdot \rangle_{\hat{g}} X \end{equation}
  be the projection operator from tangent bundle to horizontal subbundle (orthogonal to $X$).    
  
  Then  $\mathcal{L}_{X}\bar{g}$ and $A$ can be decomposed into various components: \begin{equation} \label{3.3}
 \begin{split}
 & h(\cdot,\cdot)\triangleq (\mathcal{L}_{X}\bar{g})(P\cdot, P\cdot),  \ \ \ \ \ \ \ \Lambda(\cdot,\cdot)\triangleq -u^{-2} A(P\cdot, P\cdot), \\
 & \tau(\cdot)\triangleq -u^{-2} (\mathcal{L}_{X}\bar{g})(X, P\cdot),  \ \ \ -u^{-2}A(X, P\cdot)=\tau(\cdot)-2 \nabla \log u,\end{split}
 \end{equation}
 where $h, \Lambda, \tau, \nabla\log u=d(\log u)(P(\cdot))$ are horizontal tensor fields on $\bar{M}$. By using (\ref{3.1})(\ref{3.3}) and direct computations, we obtain  
\begin{equation} \label{3.4}
\begin{split} 
D_{Y}^{\bar{g}}X&= Y(\log u)X+\langle X,Y\rangle_{\bar{g}} (\tau-\nabla \log u)-\frac{1}{2}u^2 i_{Y} \Lambda+\frac{1}{2} i_Y h\\D_{Y}^{\hat{g}}X&=Y(\log u)X-\langle X,Y\rangle_{\bar{g}} (\tau-\nabla \log u)+\frac{1}{2}u^2 i_{Y} \Lambda+\frac{1}{2} i_Y h
\end{split}
\end{equation} 
where $i_{Y} \Lambda=\Lambda_{\alpha\beta} Y^{\alpha}\hat{g}^{\beta\gamma}\frac{\partial}{\partial z^{\gamma}}$ and $i_{Y} h=h_{\alpha\beta} Y^{\alpha}\hat{g}^{\beta\gamma}\frac{\partial}{\partial z^{\gamma}}$. Here and in the followings, 
 we abuse the notations a little on  vectors and co-vectors, their meanings  should be clear from the context.

From (\ref{3.4}), we have 
 \begin{equation} \label{3.5}
 \begin{split}
 \bar{A}=&\frac{1}{2}(h-u^2 \Lambda)+X_{\ast}\otimes((\tau-\nabla \log u)-u^{-2}X(\log u)X_{\ast})+\nabla \log u \otimes X_{\ast}. \end{split}\end{equation} 
 
Using $D_{X}^{\bar{g}}Y=D_{Y}^{\bar{g}}X+[X,Y]$ and the same formula for $D^{\hat{g}}$, we obtain 
\begin{equation} \label{3.6}
\begin{split} 
D_{X}^{\bar{g}}Y=& P[X,Y]+\langle X,Y\rangle_{\bar{g}} (\tau-\nabla \log u)-\frac{1}{2}u^2 i_{Y} \Lambda+\frac{1}{2} i_Y h\\
&+[Y(\log u)+2\langle X,Y\rangle_{\bar{g}} X(\log u)u^{-2}-\tau(Y)-X\langle X, Y\rangle_{\bar{g}} u^{-2}] X\\
D_{X}^{\hat{g}}Y=& P [X,Y]-\langle X,Y\rangle_{\bar{g}} (\tau-\nabla \log u)+\frac{1}{2}u^2 i_{Y} \Lambda+\frac{1}{2} i_Y h\\
&+[Y(\log u)+2\langle X,Y\rangle_{\bar{g}} X(\log u)u^{-2}-\tau(Y)-X\langle X, Y\rangle_{\bar{g}} u^{-2}] X.
\end{split}
\end{equation}

 If $Y$ and $Z$ are horizontal, we have 
 \begin{equation} \label{3.7}
 \begin{split}
  D^{\bar{g}}_{Y}Z=& P[D^{\bar{g}}_{Y}Z]+\frac{1}{2}(u^{-2}h(Y,Z)-\Lambda(Y,Z)) X\\
  D^{\hat{g}}_{Y}Z=& P[D^{\hat{g}}_{Y}Z]+\frac{1}{2}(-u^{-2}h(Y,Z)-\Lambda(Y,Z)) X.
   \end{split}
 \end{equation}  
  Note that   $P(D^{\bar{g}}_{Y}Z)=P(D^{\hat{g}}_{Y}Z)$ if $Y$ and $Z$ are horizontal. This can be verified  from the well-known formula in Riemannian geometry: 
  \begin{equation} \nonumber
 \begin{split}
 \langle D_{U}V, W\rangle=&\frac{1}{2}(U \langle V, W\rangle+V \langle U, W\rangle-W \langle U, V\rangle\\&+ \langle [U,V], W\rangle- \langle [U,W], V\rangle- \langle [V,W], U\rangle).
 \end{split}
 \end{equation} 
 
 The difference of the connections of $\hat{g}$ and $\bar{g}$ can be summarized as follows:
  \begin{equation} \label{3.8}
\begin{split} 
D_{U}^{\hat{g}}V-D_{U}^{\bar{g}}V=& -u^{-2}h(U,V)X-\langle U,X\rangle_{\bar{g}} i_{V} \Lambda-\langle V,X\rangle_{\bar{g}} i_{U} \Lambda\\& +2\langle U,X\rangle_{\bar{g}}\langle V,X\rangle_{\bar{g}} u^{-2}(\tau-\nabla \log u) \\
\triangleq & S(U,V)
\end{split}
\end{equation}  
for all vector fields $U$ and $V$.  

The following propositions  are important for later applications. Though their  computations are  tedious,   most of them are straightforward and direct.
\begin{prop}  \label{p3.1} The Ricci curvature tensors  of $\hat{g}$  and $\bar{g}$ are related by  the following formulas: 
  \begin{equation} \label{3.9}
\begin{split}
 \hat{R}ic(X,X) =& -\bar{R}ic(X,X)+  \frac{u^4}{2}| \Lambda|^2-\frac{1}{2} |h|^2\\
&+X(\log u) tr(h)- X(tr(h))
\end{split}
\end{equation}
 \begin{equation} \label{3.10}
\begin{split}
\hat{R}ic{(X,Y)}=& -\bar{R}ic{(X,Y)}+div_{\bar{g}}(h)(Y)-Y(tr(h))\\& +tr(h) \langle \nabla\log u, Y\rangle+ h(Y,\tau-2\nabla \log u)
\end{split}
\end{equation}
 \begin{equation} \label{3.11}
\begin{split}
 \hat{R}ic({Y,Z})& =\bar{R}ic({Y,Z})-{u^2}\langle i_{Y}\Lambda, i_{Z}\Lambda \rangle -\frac{1}{2}\langle i_{Y}h, i_{Z}\Lambda\rangle-\frac{1}{2}\langle i_{Z}h, i_{Y}\Lambda\rangle 
\\& \ -u^{-2} (D_{X}^{\hat{g}} h)(Y,Z) +u^{-2} (X(\log u)-\frac{1}{2}tr(h)) h(Y,Z)  \end{split}
\end{equation}   
where $Y, Z$ are horizontal vector fields, and 
\begin{equation} \label{3.12}
div_{\bar{g}}(h)(Y)\triangleq \bar{g}^{\alpha\beta}(D^{\bar{g}}_{\alpha} h)(\partial_{\beta},Y)=\hat{g}^{\alpha\beta}(D^{\hat{g}}_{\alpha} h)(\partial_{\beta},Y).\end{equation}  \end{prop}

  \begin{prop} \label{p3.2} The Laplacian of $\log u$ w.r.t. $\hat{g}$ is given by: 
  \begin{equation} \label{3.13}
\begin{split} \triangle_{\hat{g}} \log u=& {u}^{-2} XX \log u- u^{-2} X(\log u) X(\log u)\\& +u^{-2}\bar{R}ic(X,X)-\frac{1}{4}u^{2}|\Lambda|^2
+div_{\hat{g}}(\tau)\\& +\frac{1}{4} u^{-2}|h|^2+\frac{1}{2} u^{-2} X(tr (h)).
\end{split}
\end{equation}
\end{prop} 

\begin{prop} \label{p3.3}The 2-form $\Lambda$ satisfies:
  \begin{equation} \label{3.14}
\begin{split} d\Lambda &= -d\tau\wedge u^{-2}{X}_{\ast}-\tau\wedge \Lambda\\
 div_{\hat{g}}(u^2\Lambda)(Y)&=2\bar{R}ic(X,Y)+ Y(tr (h))-h(\tau-2\nabla \log u, Y)\\ & \ \ \ +u^2\Lambda(\tau, Y)- tr(h) \langle \nabla \log u,Y\rangle-div_{\hat{g}}(h)(Y),
\end{split}
\end{equation}
where ${X}_{\ast}$ is the co-vector field dual to $X$ w.r.t. $\bar{g}$.
\end{prop}

 To prove Proposition \ref{p3.1},    we need to  compute the curvature tensor: \begin{equation}
 \begin{split} \label{3.15}
 \bar{R}(U,V,W, Z)&\triangleq -\langle D^{\bar{g}}_UD^{\bar{g}}_VW-D^{\bar{g}}_VD^{\bar{g}}_UW-D^{\bar{g}}_{[U,V]}W, Z \rangle_{\bar{g}}.\end{split}
 \end{equation}
 Let $W=X$,  we obtain 
 \begin{equation} \label{3.16}
  \begin{split}
 \bar{R}(U,V,X, Z)&=-(D^{\bar{g}}_U\bar{A})(V,Z)+(D^{\bar{g}}_V\bar{A})(U,Z).\end{split}
 \end{equation}

    Taking trace on $V$ and $Z$ in (\ref{3.16}), we find  
    \begin{equation} \label{3.17}
  \begin{split}
 \bar{R}ic(X, U)=-U tr_{\bar{g}}(\bar{A})+(D^{\bar{g}}_{\alpha} \bar{A})(U,\partial_{\beta})\bar{g}^{\alpha\beta}.
 \end{split}
 \end{equation} 
 Let $U=X$ in (\ref{3.17}),  we have 
 \begin{equation} \label{3.18}
  \begin{split}
 \bar{R}ic(X, X)=-X tr_{\bar{g}}(\bar{A})-tr_{\bar{g}}(\bar{A}^2)+div_{\bar{g}}(\bar{A}(X, \cdot)).
 \end{split}
 \end{equation}  
 
 From (\ref{3.5}), we have  \begin{equation} \label{3.19}
  \begin{split}
 \bar{A}(X, \cdot)=\langle -u^2(\tau-\nabla\log u) +X(\log u)X, \cdot \rangle\end{split}
 \end{equation} 
 and 
 \begin{equation} \label{3.20}
  \begin{split} tr(\bar{A})&=\frac{1}{2} tr(h)+X(\log u)\\
  tr_{\bar{g}}(\bar{A}^2)&=\frac{1}{4}|h|^2-\frac{u^4}{4}|\Lambda|^2-2u^2\langle \tau-\nabla\log u, \nabla \log u\rangle+(X\log u)^2.  \end{split}
 \end{equation}   
 Substituting (\ref{3.19})  and (\ref{3.20}) into (\ref{3.18}), we obtain 
 \begin{equation} \label{3.21}
 \begin{split}
 \bar{R}ic(X,X)=&u^2 div_{\bar{g}}(\nabla \log u-\tau)-\frac{1}{4}|h|^2+\frac{u^4}{4}|\Lambda|^2\\ &+\frac{1}{2} X(\log u)tr(h)-\frac{1}{2} X(tr (h)). \end{split}\end{equation}

 Let $Y$ be horizontal, $E_i$, $i=1,2,3$ be a horizontal orthonormal basis, using (\ref{3.17}) and (\ref{3.5}),  we obtain
  \begin{equation} \label{3.22}
 \begin{split}
 \bar{R}ic(X,Y)=& -\frac{1}{2} Y(tr (h))-YX(\log u)+\frac{1}{2} div_{\bar{g}}(h+u^2\Lambda)(Y)\\
 &+\langle D^{\bar{g}}_{E_i}X,Y\rangle \langle\tau-\nabla \log u,E_i\rangle-u^{-2} \langle D^{\bar{g}}_{X}X,Y\rangle_{\bar{g}}X(\log u)\\
 &+\langle D^{\bar{g}}_{X} \nabla \log u,Y\rangle_{\bar{g}} +\langle \nabla \log u,Y\rangle (\frac{1}{2}tr(h)+X(\log u))\\
  =&\frac{1}{2} div_{\bar{g}}(h+u^2\Lambda)(Y)-\frac{1}{2} Y(tr (h))+\frac{1}{2}h(\tau-2\nabla \log u, Y)\\ & -\frac{1}{2}u^2\Lambda(\tau, Y)+\frac{1}{2} tr(h) \langle \nabla \log u,Y\rangle. \end{split}\end{equation}  
  
  Now we  compute the horizontal components $\bar{R}ic(Y,Z)$, where $Y,Z$ are horizontal vector fields.  From (\ref{3.16}), we have 
  \begin{equation} \label{3.23}
  \begin{split}
  \bar{R}(X,Y,X,Z)&=-(D^{\bar{g}}_X\bar{A})(Y,Z)+(D^{\bar{g}}_Y\bar{A})(X,Z) \\
  =&-\frac{1}{2}D^{\bar{g}}_{X}(h-u^2\Lambda)(Y,Z)-u^2 \langle D^{\hat{g}}_{Y}(\tau-\nabla \log u), Z\rangle\\
  &-\frac{1}{4}\langle i_{Y}(h-u^2\Lambda), i_{Z}(h+u^2\Lambda)\rangle+\frac{1}{2}X(\log u)(h-u^2\Lambda)(Y,Z) \\
  &+u^2\langle \tau-\nabla \log u, Y\rangle\langle \tau-\nabla \log u, Z\rangle.   \end{split}\end{equation}
  Using the symmetry on $Y$ and $Z$ in $\bar{R}(X,Y,X,Z)$, we obtain 
   \begin{equation} \label{3.24}
  \begin{split}
  \bar{R}(X,Y,X,Z)=&-\frac{1}{2}u^2 \langle D^{\hat{g}}_{Y}(\tau-\nabla \log u), Z\rangle-\frac{1}{2}u^{2}\langle D^{\hat{g}}_{Z}(\tau-\nabla \log u), Y\rangle\\
  &-\frac{1}{4}\langle i_{Y}h, i_{Z}h\rangle+\frac{1}{4}u^4\langle i_{Y}\Lambda, i_{Z}\Lambda\rangle +\frac{1}{2}X(\log u)h(Y,Z) \\
  &+u^2\langle \tau-\nabla \log u, Y\rangle\langle \tau-\nabla \log u, Z\rangle-\frac{1}{2}(D^{\bar{g}}_{X}h)(Y,Z).  
   \end{split}\end{equation}  
  If $U,V,Y,Z$ are horizontal, using (\ref{3.15}),   we have 
 \begin{equation} \label{3.25}
  \begin{split}
 \bar{R}(Y,U, Z,V)=&R(Y,U,Z,V)+[\bar{A}(Y,Z)\bar{A}(U,V)-\bar{A}(Y,V)\bar{A}(U,Z)]u^{-2}\\
 &-\Lambda(Y,U)\bar{A}(Z,V)\end{split}
 \end{equation} 
 where 
 \begin{equation} \label{3.26}
 \begin{split}{R}(Y,U,Z,V)=&-\langle P D^{\bar{g}}_Y(P(D^{\bar{g}}_UZ))-
  P D^{\bar{g}}_U(P(D^{\bar{g}}_YZ))-P(D^{\bar{g}}_{P[Y,U]}Z), V\rangle_{\bar{g}}\\
  &-\Lambda(Y,U)\langle P[X,Z],V\rangle_{\bar{g}}
    \end{split}\end{equation}  
    is a horizontal $(0,4)$-tensor field. 
  Using (\ref{3.24}) and (\ref{3.25}), we  have 
  \begin{equation} \label{3.27}
 \begin{split}\bar{R}ic(Y,Z)=&\bar{g}^{\alpha\beta} \bar{R}(Y,
 P\partial_{\alpha},Z,P\partial_{\beta})-u^{-2}  \bar{R}(X,Y,X,Z)\\
 =& \bar{g}^{\alpha\beta} {R}(Y, 
 P\partial_{\alpha},Z,P\partial_{\beta}) +\frac{1}{4} tr(h)(u^{-2}h-\Lambda)(Y,Z)\\
 &+\frac{1}{2} \langle D^{\hat{g}}_{Y}(\tau-\nabla \log u), Z\rangle+\frac{1}{2}\langle D^{\hat{g}}_{Z}(\tau-\nabla \log u), Y\rangle\\
  &+\frac{1}{2}u^2\langle i_{Y}\Lambda, i_{Z}\Lambda\rangle -\frac{1}{2}u^{-2}X(\log u)h(Y,Z) \\
  &-\langle \tau-\nabla \log u, Y\rangle\langle \tau-\nabla \log u, Z\rangle+\frac{1}{2}u^{-2}(D^{\bar{g}}_{X}h)(Y,Z)\\
  &-\frac{1}{4}  \langle i_{Y}h, i_{Z}\Lambda\rangle -\frac{1}{4}  \langle i_{Z}h, i_{Y}\Lambda\rangle. \end{split}\end{equation}

  \begin{rem} \label{r3.1}
 Same computations can be carried out for the Riemannian metric $\hat{g}$ defined in (\ref{1.2}). It turns out that the corresponding formulas of $\hat{g}$ for the connection coefficients,  full curvature tensor and Ricci curvature tensor 
 can be  obtained by simply substituting $u$ with $\sqrt{-1}u$ in  the formulas of $\bar{g}$. 
  \end{rem}
 \begin{pf}  of Proposition \ref{p3.1}.\\
  The results   follow from  (\ref{3.21}) (\ref{3.22}) (\ref{3.27}) and Remark \ref{r3.1}
 \end{pf}
 \begin{pf} of Proposition \ref{p3.2}. \\
 Since $d\log u=\nabla \log u-X(\log u) u^{-2} X_{\ast}$, we have \begin{equation} \label{3.28}
 \triangle_{\hat{g}} \log u=div_{\hat{g}}(\nabla \log u)-div_{\hat{g}}(X(\log u)u^{-2}X_{\ast}).
 \end{equation}
 By direct computations, we have 
 \begin{equation} \label{3.29}
 \begin{split}
 div_{\hat{g}}(-X(\log u) u^{-2}X_{\ast}) =&XX(\log u)u^{-2}-X(\log u)X(\log u) u^{-2}\\&+\frac{1}{2} tr(h) X(\log u)u^{-2}.
 \end{split}
 \end{equation}
 
 Combining (\ref{3.28}) (\ref{3.29}) and  (\ref{3.21}), Proposition \ref{p3.2} follows.
 \end{pf}
 \begin{pf} of Proposition \ref{p3.3}. 
 
  Let $\{E_1, E_2, E_3, E_0=u^{-1}X\}$ be a  normal basis of $\hat{g}$, $\{E^{\ast}_{\alpha}\}$ be  its dual, where $E_0=-u^{-1}X_{\ast}$.  $\Lambda$ is the  horizontal two form 
 $$
 \Lambda=\sum_{i<j} \Lambda(E_i,E_j) E^{\ast}_i\wedge E^{\ast}_j.$$
 
  From (\ref{3.2})(\ref{3.3}),    we have 
 \begin{equation} \label{3.30}
 d(u^{-2} X_{\ast})=-\Lambda-\tau\wedge u^{-2}X_{\ast}.
 \end{equation}
 Taking exterior differential on (\ref{3.30}), we get
  \begin{equation} \label{3.31}
 d\Lambda+d\tau\wedge u^{-2}X_{\ast}+\tau\wedge \Lambda=0,
 \end{equation} 
  which is the first equation in (\ref{3.14}). The second equation in (\ref{3.14}) can be derived from (\ref{3.22}).  
 

\end{pf}

\section{Groupoid approach}  
     
   In some  sense,  we  need certain smooth enough  limit space on which the Einstein equations can be written down. In general, because of the collapsing, the Gromov-Hausdorff convergence in Section 2 can not be improved to usual $C^{\alpha}$ convergence in Cheeger-Gromov sense. 
 
 Our  strategy is  firstly to consider the problem  in each  fixed  coordinate patch, and prove that there is  a  subsequence of $\bar{g}_j$, $\hat{g}_j$ and $X_j$ which convergences to a smooth limit $\bar{g}_{\infty}, \hat{g}_{\infty}, X_{\infty}$ in  $C^{\alpha}$ and weak $W^{1,p}$ topologies.  Then we proceed to investigate how the  limits  in these coordinate patches can be glued together, with an aim to  prove that all these  limits  are actually flat.

 \subsection{$C^{\alpha}$-convergence}

By (\ref{2.9}),  $r^{w,h}(W^{1,p_0}, \epsilon_0)(\bar{x}_j)=1$,   there is a smooth  immersion  $\Psi$ defined on a neighborhood of $\{z\in \mathbb{R}^4: |z|\leq 1/3\}$ with $\Psi(0)=P$ such that (\ref{2.1}) (\ref{2.2}) and (\ref{2.3}) hold for $r=1/3$. Now we pull back everything  to  $\Omega_0=\{ |z|\leq 1/3\}$, and use the same notations to denote them, e.g.  $\hat{g}\triangleq \Psi^{\ast} \hat{g}_j$, $\bar{g}_j\triangleq \Psi^{\ast} \bar{g}_j$, $X_j\triangleq \Psi^{\ast}X_j$.   The vector field $X_j$ can also be pulled back  since the map  $\Psi$ is an immersion. We consider the convergence of the  quantities $(\hat{g}_j, \bar{g}_j, X_j)$ on $\Omega_0$. 

Using $\partial \log u_j^2=-2\langle D^{\bar{g}_j} X_j, X_j\rangle_{\bar{g}_j}/ u_j^{2}$
and (\ref{2.3}), we know $||\partial \log u_j||_{L^{p_0}(\Omega_0)} \leq C$.  Since   $u_j(0)=1$ and (\ref{2.2}) (\ref{2.3}),   we  have  $|| \log u_j||_{W^{1,p_0}(\Omega_0)}$, $|| X_j||_{W^{1,p_0}(\Omega_0)}$,  $||\hat{g}_j||_{W^{1,p_0}(\Omega_0)}$ and $||\bar{g}_j||_{W^{1,p_0}(\Omega_0)}$  are uniformly bounded.  By Sobolev imbedding theorem, for any $0<\alpha< \alpha_0=1-\frac{n}{p_0}$, there is a subsequence of $j$ so that $\hat{g}_j$, $\bar{g}_j$, $X_j$, $u_j$ converge to $\hat{g}_{\infty}$, $\bar{g}_{\infty}$, $X_{\infty}$, $u_{\infty}$ in $C^{\alpha}(\Omega_0)$ topology. Since the space $W^{1,p_0}(\Omega_0)$ is reflexive, we may assume the convergences also take place  in weak $W^{1,p_0}(\Omega_0)$ topology.  

In particular, multiplying the equation  $\triangle_{\hat{g}_j}z^{\gamma}=0$ with a $L^{p^{\ast}_0}$ function ($\frac{1}{p_0}+\frac{1}{p^{\ast}_0}=1$), integrating  and taking limits,  we have 
 \begin{equation} \label{4.1}\triangle_{\hat{g}_{\infty}}z^{\gamma}=-\hat{g}_{\infty}^{\alpha\beta} \Gamma(\hat{g}_{\infty})^{\gamma}_{\alpha\beta}=0
\end{equation}
 as $L^{p_0}$ functions, i.e. $\{z^{\beta}\}$ is still a harmonic coordinate system for $\hat{g}_{\infty}$.   
 
  \begin{thm}\label{t4.1}
  $\hat{g}_{\infty}$, $\bar{g}_{\infty}$, $X_{\infty}$, $u_{\infty}$ are smooth on $\{|z|<1/3\}$, and 
  $$\mathcal{L}_{X_{\infty}}\bar{g}_{\infty}=\mathcal{L}_{X_{\infty}}\hat{g}_{\infty}=\mathcal{L}_{X_{\infty}}{u}_{\infty}=\mathcal{L}_{X_{\infty}}{\Lambda}_{\infty}=0.$$
  
  Moreover, for any  $k\in \mathbb{Z}_{+}$ and $0<\delta<3$, there is a constant $C_{k,\delta}>0$ depending only on $k, \delta$ such that 
  \begin{equation} \label{4.2}
  ||\hat{g}_{\infty},\bar{g}_{\infty}, X_{\infty}, u_{\infty},{\Lambda}_{\infty} ||_{C^{k}(|z|\leq 3 -\delta)}\leq C_{k,\delta}.
  \end{equation}
  \end{thm}

\begin{pf}
Multiplying (\ref{3.13}) with a smooth function $\xi$ with compact support, integrating by parts (with volume form of $\hat{g}_j$), and taking limits, we have  
\begin{equation} \label{4.3}
\int (\hat{g})_{\infty}^{\alpha\beta} (\log u_{\infty})_{\alpha} \xi_{\beta} dv_{\hat{g}_{\infty}}= \int f_{\infty} \xi dv_{\hat{g}_{\infty}}, 
\end{equation}
  for some function $f_{\infty}\in L^{\frac{p_0}{2}}$. Here we can not say $f_{\infty}$ is  $\frac{1}{4}u_{\infty}^{4}|\Lambda_{\infty}|^2$,  $f_{\infty}$  is only  the weak limit of $\frac{1}{4}u_j^{4}|\Lambda_j|^2$ in $L^{\frac{p_0}{2}}$.   $\Lambda_j$ is weakly converging  to $\Lambda_{\infty}$ in $L^{{p_0}}$, it does not imply that $ |\Lambda_j|^2$ is  weakly converging  to $ |\Lambda_{\infty}|^2$ in $L^{\frac{p_0}{2}}$. The same problem also appears in the following equations  of $\hat{g}_{\infty}$ and $\Lambda_{\infty}$. 

That is to say, $\triangle_{\hat{g}_{\infty}} \log u_{\infty}=-f_{\infty}\in L^{\frac{p_0}{2}}$ in weak sense. 

In harmonic coordinates, the Laplace operator takes the form 
$\triangle_{\hat{g}_{\infty}} = (\hat{g})_{\infty}^{\alpha\beta}\partial_{\alpha} \partial_{\beta}$. Since  $(\hat{g})_{\infty}^{\alpha\beta}\in C^{\alpha} \cap W^{1,p_0}$, $f_{\infty}\in L^{\frac{p_0}{2}}$, $u_{\infty}\in W^{1,p_0}$,  by $L^{p}$ estimate of  elliptic equations (see \cite{GT}), we obtain  $u_{\infty} \in W_{loc}^{2,\frac{p_0}{2}}$. 
Now we shall apply  the same argument to the equation of $\hat{g}_j$ which is  the  Ric curvature $Ric(\hat{g}_j)$ of $\hat{g}_j$:   
\begin{equation} \label{4.4}
-2 Ric(\hat{g}_j)_{\alpha\beta}=\hat{g}^{\gamma\delta}_j \partial_{\gamma}\partial_{\delta}(\hat{g}_j)_{\alpha\beta}+(\hat{g}^{-2}_j\ast \partial \hat{g}_j \ast \partial \hat{g}_j)_{\alpha\beta} 
\end{equation}
in harmonic coordinates $\{z^{\gamma}\}$.

 To simplify notations, we omit the subscript $j$. Using $\partial_{\alpha}=P\partial_{\alpha}+\langle \partial_{\alpha}, X\rangle_{\hat{g}}u^{-2}X$, 
where $P\partial_{\alpha}$ is the horizontal projection of $\partial_{\alpha}$,  substituting (\ref{3.9}) (\ref{3.10}) and (\ref{3.11})  into  $Ric(\hat{g})_{\alpha\beta}$, we have 
\begin{equation} \label{4.5}
\begin{split}
 Ric(\hat{g})_{\alpha\beta}&=I+II+III
  \end{split}
\end{equation}    
where  
\begin{equation} \label{4.6}
\begin{split} I=&Ric(\bar{g})(P\partial_{\alpha},P\partial_{\beta})-u^{-4}X_{\alpha}X_{\beta} Ric(\bar{g})(X,X)\\
 & -u^{-2}Ric(\bar{g})(X, P\partial_{\alpha})X_{\beta}-u^{-2}Ric(\bar{g})(X, P\partial_{\beta})X_{\alpha}\end{split}
\end{equation} 
   \begin{equation}
\begin{split} \label{4.7}
II=&u^{-4}X_{\alpha}X_{\beta} [\frac{u^4}{2}| \Lambda|^2-\frac{1}{2} |h|^2+X(\log u) tr(h))]\\
&+u^{-2}[tr(h) \langle \nabla\log u, \partial_{\alpha}\rangle+ h(\partial_{\alpha},\tau-2\nabla \log u)]X_{\beta}\\&+u^{-2}[tr(h) \langle \nabla\log u, \partial_{\beta}\rangle+ h(\partial_{\beta},\tau-2\nabla \log u)]X_{\alpha} \\&
-{u^2}\langle i_{\partial_{\alpha}}\Lambda, i_{\partial_{\beta}}\Lambda \rangle -\frac{1}{2}\langle i_{\partial_{\alpha}}h, i_{\partial_{\beta}}\Lambda\rangle-\frac{1}{2}\langle i_{\partial_{\beta}}h, i_{\partial_{\alpha}}\Lambda\rangle 
 \\& +u^{-2} (X(\log u)-\frac{1}{2}tr(h)) h(\partial_{\alpha},\partial_{\beta}) \end{split}
\end{equation}  
\begin{equation} \label{4.8}
\begin{split} III&=u^{-4}X_{\alpha}X_{\beta} [- X(tr(h))]+u^{-2}[div_{\bar{g}}(h)(\partial_{\alpha})-(P\partial_{\alpha})(tr(h))] X_{\beta}\\&+u^{-2}[div_{\bar{g}}(h)(\partial_{\beta})-(P\partial_{\beta})(tr(h))]X_{\alpha} 
-u^{-2} (D_{X}^{\hat{g}} h)(P\partial_{\alpha},P\partial_{\beta}).
\end{split}
\end{equation}

It is not hard to see 
\begin{equation} \label{4.9}
\begin{split} 
III=&u^{-1}(\partial h+\hat{g}^{-1}\ast \partial \hat{g}\ast h)\ast \hat{g}^{-1} \ast [u^{-1}X_{\ast} +(u^{-1}X_{\ast})^{3}\ast \hat{g}^{-1}\\&+(u^{-1}X_{\ast})^{5}\ast \hat{g}^{-2}]\end{split}
\end{equation} 
and 
\begin{equation} \label{4.10}
\begin{split} 
|I| & \preceq |\bar{R}ic|\\
|II|& \preceq u^{-2}|\mathcal{L}_{X}\bar{g}|^2+u^{-1}|\mathcal{L}_{X}\bar{g}|(|\nabla \log u|+u|\Lambda|)+u^2|\Lambda|^2.  \end{split}
\end{equation}  

Here and in the followings,  $A \preceq B$ means that there is a universal constant $C$ such that $A \leq CB$.    
Multiplying (\ref{4.4}) with a smooth function $\psi$ with compact support, integrating by parts with volume form $dv_{\hat{g}}$, we have 
 \begin{equation} \label{4.11}
 \begin{split}
& |\int \hat{g}^{\gamma\delta} \partial_{\gamma}(\hat{g})_{\alpha\beta} \partial_{\delta} \psi+[(\hat{g}^{-2}\ast \partial \hat{g} \ast \partial \hat{g})_{\alpha\beta}+u\Lambda\ast u\Lambda \ast \hat{g}^{-1}] \psi dv_{\hat{g}}|\\&
\preceq \int [|\bar{R}ic|+u^{-1}|\mathcal{L}_{X}\bar{g}|(u^{-1}|\mathcal{L}_{X}\bar{g}|+|\nabla \log u|+u|\Lambda|+u^{-1}|\partial X|+|\partial \hat{g}|)] \psi\\
&+u^{-1}|h| |d \psi| dv_{\hat{g}}. \end{split}
   \end{equation}
   Taking limits, the right hand side of (\ref{4.11}) converges to zero, and there is a $L^{\frac{p_0}{2}}$ tensor field $f_{\alpha\beta}$ which is the weak limit of $\hat{g}^{-2}\ast \partial \hat{g} \ast \partial \hat{g}+u\Lambda\ast u\Lambda\ast \hat{g}^{-1}$ in $L^{\frac{p_0}{2}}$ such that 
   \begin{equation} \label{4.12}
   \hat{g}_{\infty}^{\gamma\delta} \partial_{\gamma}\partial_{\delta}(\hat{g}_{\infty})_{\alpha\beta} =f_{\alpha\beta}.
      \end{equation}
   
 Applying the standard $L^{p}$ estimate to (\ref{4.12}) yields  
\begin{equation} \label{4.13}
(\hat{g}_{\infty})_{\alpha\beta}\in W_{loc}^{2,\frac{p_0}{2}}.
\end{equation}
 The regularity of $\Lambda_{\infty}$ can be done better. Wedging the first formula of  (\ref{3.14})  with  a smooth form $\xi$ with compact support and taking limits, we find 
 \begin{equation} \label{4.14}
 \begin{split}
 -\int  \Lambda \wedge d \xi&=\int  d\tau\wedge u^{-2}X_{\ast} \wedge \xi-\tau\wedge \Lambda\wedge \xi \\
 &=-\int  \tau\wedge u^{-2}X_{\ast} \wedge d \xi \rightarrow 0.\end{split}
 \end{equation}
 Here we have used 
    \begin{equation} \label{4.15}
 d(u^{-2} X_{\ast})=\Lambda-\tau\wedge u^{-2}X_{\ast}.
 \end{equation} 
 That  is to say,  we have 
  \begin{equation} \label{4.16} d\Lambda_{\infty}=0
  \end{equation}
  in weak sense. Similar argument on the second formula of (\ref{4.14}) yields  
  \begin{equation} \label{4.17} div_{\hat{g}_{\infty}}(u_{\infty}^2\Lambda_{\infty})=0
  \end{equation}   
  in weak sense.   Since $\hat{g}_{\infty}, u_{\infty}\in W_{loc}^{2,\frac{p_0}{2}}$, the elliptic regularity on equations (\ref{4.16}) and (\ref{4.17}) imply that $\Lambda_{\infty}\in W_{loc}^{2,\frac{p_0}{2}}$.

 The equation for $X_{\infty}$ can be obtained directly by taking weak  limits in $L^{p_0}$ in (\ref{3.4}): 
  \begin{equation} \label{4.18}
\begin{split} D_{Y}^{\bar{g}_{\infty}}X_{\infty}&=Y(\log u_{\infty})X_{\infty}-\langle X_{\infty},Y\rangle_{\bar{g}_{\infty}} \nabla \log u_{\infty}-\frac{1}{2}u_{\infty}^2 i_{Y} \Lambda_{\infty}\\D_{Y}^{\hat{g}_{\infty}}X_{\infty}&=Y(\log u_{\infty})X_{\infty}-\langle X_{\infty},Y\rangle_{\hat{g}_{\infty}} \nabla \log u_{\infty}+\frac{1}{2}u_{\infty}^2 i_{Y} \Lambda_{\infty}
\end{split}
\end{equation}  
in weak  sense.    For instance, in  (\ref{3.4}),  $Y(\log u)$ weakly converges  to $Y(\log u_{\infty})$ in $L^{p_0}$, $X$ converges  to $X_{\infty}$ in $C^{\alpha}$, hence $Y(\log u)X$ weakly converges  to $Y(\log u_{\infty}) X_{\infty}$ in $L^{p_0}$.  From  (\ref{4.18}), we know $X_{\infty}\in W_{loc}^{2,\frac{p_0}{2}}$. 

 In the next step, we shall improve the regularities of $\log u_{\infty}$, $\bar{g}_{\infty}$, $\hat{g}_{\infty}$.   The point is that we can  recalculate their equations,  instead of using the equations before taking limits.

Note that $X_{\infty} \in W_{loc}^{2,\frac{p_0}{2}}$, $\mathcal{L}_{X_{\infty}}\bar{g}_{\infty}=\mathcal{L}_{X_{\infty}}\hat{g}_{\infty}=0$ (in $W_{loc}^{1,\frac{p_0}{2}}$). Since $\bar{g}_{\infty}, \hat{g}_{\infty} \in W_{loc}^{2, \frac{p_0}{2}}$,    one can calculate the Ricci curvatures of $\hat{g}_{\infty}$, $\hat{g}_{\infty}$, and equation of $\log u_{\infty}$: 
\begin{equation} \label{4.19}
\begin{split}
Ric(\bar{g}_{\infty})_{\alpha\beta}&=0\\
Ric(\hat{g}_{\infty})_{\alpha\beta}&=  \frac{1}{2} |\Lambda_{\infty}|^2(X_{\infty})_{\alpha}(X_{\infty})_{\beta}
-{u_{\infty}^2}\langle i_{\partial_{\alpha}}\Lambda_{\infty}, i_{\partial_{\beta}}\Lambda_{\infty} \rangle \\
\triangle_{\hat{g}_{\infty}} \log u_{\infty}&=-\frac{1}{4}u_{\infty}^{2}|\Lambda_{\infty}|^2
\end{split}
\end{equation}
which hold (in $L_{loc}^{\frac{p_0}{2}}$) in the sense of integration by parts. Here the reason why one can identity the righthand side of (\ref{4.19}) is because of (\ref{4.18}). 
Now a bootstrap argument based on (\ref{4.19})(\ref{4.16})(\ref{4.17})(\ref{4.18}) implies  that all these quantities $\hat{g}_{\infty}$, $\bar{g}_{\infty}$, $\Lambda_{\infty}$, $u_{\infty}$, $X_{\infty}$ are smooth in $\hat{g}_{\infty}$-harmonic coordinates  $\{z^{\alpha}\}$.

The estimate (\ref{4.2}) holds because all the above elliptic estimates are uniform on  any  compact sets  of the domain $\{|z|<\frac{1}{3}\}$. 
This completes the proof.    
\end{pf}

 \subsection{Strong maximum principle}
 
 Let $(L, d_{L}, \bar{x})$ be the limit metric space obtained in Section 2.   After passing to a subsequence, we may assume the sequence $(\bar{M}_j, \hat{g}_j, \bar{x}_j)$ is  chosen so that 
for any $j\in \mathbb{N}$,  there are  approximating  maps $I^{j}_{k}: B_{L}(\bar{x}, j)\rightarrow B_{\hat{g}_k}(\bar{x}_k, j+\frac{1}{k})\subset \bar{M}_k$ and $J^{j}_{k}: B_{\hat{g}_k}(\bar{x}_k, j) \rightarrow  B_{L}(\bar{x}, j+\frac{1}{k}) $ for any $k\geq j$,   such that $I^{j}_{k}(\bar{x})=\bar{x}_k$, $J^{j}_{k}(\bar{x}_k)=\bar{x}$ and the following two conditions hold: 

i)  the $2/k$-neighborhoods of $I^{j}_{k}(B_{L}(\bar{x}, j))$ and $J^{j}_{k}(B_{\hat{g}_k}(\bar{x}_k, j)$ contain $B_{\hat{g}_k}(\bar{x}_k, j+\frac{1}{k})$ and $B_{L}(\bar{x}, j+\frac{1}{k})$ respectively; 

ii)  for any $x,y\in B_{L}(\bar{x}, j)$, $\tilde{x},\tilde{y}\in B_{\hat{g}_k}(\bar{x}_k, j)$,  we have 
\begin{equation}
\begin{split}
& |d_{L}(x,y)-d_{\hat{g}_k}(I^{j}_{k}(x), I^{j}_{k}(y)| <1/k\\
& |d(\tilde{x},\tilde{y})-d_{L}(J^{j}_{k}(\tilde{x}), J^{j}_{k}(\tilde{y})| <1/k\\
& d_{L}(x,J^j_{k}I^j_{k}(x))+d_{\hat{g}_k}(\tilde{x}, I^{j}_{k}J^{j}_{k}(\tilde{x}))<1/k.
\end{split}
\end{equation}

On the other hand,  there are  an  increasing sequence 
$n_1\leq n_2\leq n_3\cdots$ of natural numbers and a sequence of points $x_0=\bar{x}, x_1, x_2 \cdots $ in $L$ such that 

i) $d_{L}(x_i,x_j) \geq 1/10$, for $i\neq j$; 

ii) for each $l\in \mathbb{N}$,   $\{x_i\}_{i=0}^{n_l}$ is a set of maximal number of points in $B_{L}(\bar{x}, l)$ with mutual distances $\geq 1/10$. 

From Section 2, $n_l\leq Ce^{Cl}$ for some universal constant $C$. 

It is clear that $\{B_{L}(x_i, 1/20)\}_{i=0}^{\infty}$ are disjoint, and 
$\cup_{i=0}^{n_l} B_{L}(x_i, 1/10)\supset B_{L}(x_i, l) $,  for each $l\in \mathbb{N}$. 

For each $k\in \mathbb{N}$, let $x^k_{i} \triangleq I^{k}_k(x_i) \in \bar{M}_k$, for $0\leq i\leq  n_k$. 
Then $\{B_{\hat{g}_k}(x^{k}_i, 1/9)\}_{i=0}^{n_k}$ is a covering of $B_{\hat{g}_k}(\bar{x}_{k}, k)$. Since $r^{w,h}(W^{1,p_0}, \epsilon_0)(x^{k}_i) \geq 1/2$, for all $0\leq i \leq n_k$ (we assume $m_k>> k$ by choosing the subsequence), there exist a family of immersions $\Psi^{k}_i: \{ |z|<1/2\} \rightarrow \bar{M}_k$  with $\Psi^{k}_i(0)=x^{k}_i$ such that (\ref{2.1}) (\ref{2.2})(\ref{2.3}) hold for $r=1/2$, for all $0\leq i\leq n_k$. 

For fixed $i$, since the pullbacks of $\hat{g}_{k}$, $\bar{g}_k$ and $X_k$ using  $\Psi^{k}_i$ have uniformly bounded   $W^{1,p_0}$ norms on $\{ |z|<1/2\} \triangleq \Omega_i$,   by Sobolev imbedding theorem, there is a subsequence of   $\hat{g}_{k}$, $\bar{g}_k$ and $X_k$ which converges  to $\hat{g}_{\infty}$, $\bar{g}_{\infty}$ and $X_{\infty}$ in $C^{\alpha}$ and weak $W^{1,p_0}$ topologies. 

Using a diagonal process, one may assume the subsequence is chosen so that it converges on any $\Omega_i$ for all $0\leq i <\infty$. 

Let $\bar{M}_{\infty}= \amalg_{i=0}^{\infty} \Omega_i$ be the manifold of disjoint union of these $\Omega_i$. By Theorem \ref{t4.1},  $\hat{g}_{\infty}$, $\bar{g}_{\infty}$ and $X_{\infty}$ are smooth on $\bar{M}_{\infty}$. 

\begin{thm} \label{t4.2}
Both $\hat{g}_{\infty}$ and $\bar{g}_{\infty}$ are flat, $X_{\infty}$ is parallel with respect to both $\hat{g}_{\infty}$ and  $\bar{g}_{\infty}$ on $\bar{M}_{\infty}$.
\end{thm}

The idea is  to consider the convergence  of those  gluing maps between these $\Omega_i$.   It turns out that the appropriate framework for this   discussion is  to  consider groupoids,  instead of manifolds.    

\begin{pf}
 Note that  the horizontal subbundle  has dimension 3,  if one arbitrarily fix a (local) orientation,  the  horizontal Hodge dual $\omega$ of $u^3\Lambda$ can be defined:
 \begin{equation} \label{4.21}
 \omega=u^3 \ast_{3}
 \Lambda. \end{equation}
   Let  $\omega_{\infty}=u_{\infty}^3\ast_{3}\Lambda_{\infty}$  be  the  corresponding one-form   on the limit.   
  
For any $0<\delta\leq 10^{-1}$,  we assert  that there is a universal constant $C>0$ such that 
\begin{equation} \label{4.22}
\sup_{i}\sup_{\{|z|\leq (1-\delta) 1/2, z\in \Omega_i\}} 4 | \nabla \log u_{\infty}|^2+ |\omega_{\infty}|^2 u_{\infty}^{-4} \leq C \delta^{-2}.
\end{equation}

 For any point $z\in \Omega_i$ with $|z|\leq (1-\delta) 1/2$, we have $B_{\hat{g}_{\infty}}(z, C^{-1} \delta)\subset \Omega_i$ for some universal constant $C$.  The interior estimate of stationary solutions (see \cite{A1},  \cite{Chen1} Theorem 1.3) says  that  
$$ \sup_{B_{\hat{g}_{\infty}}(z, \frac{1}{2}C^{-1} \delta)} 4 | \nabla \log u_{\infty}|^2+ |\omega_{\infty}|^2 u_{\infty}^{-4} \leq C \delta^{-2},$$  which particularly gives the estimate at $z$. (\ref{4.22}) is proved.

Fix a $\delta=\delta_0\leq 10^{-1}$, let $\Omega^{\prime}_i=\{|z|<1/2(1-\delta_0):z\in \Omega_i\}$, $M^{\prime}=\cup_{i}\Omega^{\prime}_i\subset \bar{M}_{\infty}$.  Then  $$f(x)\triangleq 4 | \nabla \log u_{\infty}|^2(x)+ |\omega_{\infty}|^2(x)$$ is a smooth and bounded function on $M^{\prime}$. Denote $\bar{f}=\sup_{x\in M^{\prime}} f(x)$. 

We claim $\bar{f}=0$.

Once we have proved $f(x) \equiv 0$, then it is easy to see that $\hat{g}_{\infty}$, $\bar{g}_{\infty}$ are flat, and $X_{\infty}$ are parallel. 

Suppose the claim is not true, i.e.  $\bar{f}>0$,  let  $\{y_k\in M^{\prime}\}$ be a sequence of points with  $f(y_k)\geq (1-k^{-1}) \bar{f} \rightarrow \bar{f}$, as $k\rightarrow \infty$.  For each $y_k$, there is an $m(k)>0$ with  $y_k\in \Omega^{\prime}_{m(k)}$.

For fixed  $k=k_0$, we argue as follows. 

 For  each $i=1,2,\cdots$,  let $\{\Psi^{i}_{j}: \Omega_j \rightarrow \bar{M}_i\}^{n_i}_{j=0}$ be the family of immersions such that $\cup_{j=0}^{n_i}\Psi^i_j(|z|<1/8)$ covers  $B_{\hat{g}_i}(\bar{x}_i,i)$.  Recall that the  increasing sequence  $n_1\leq n_2 \cdots$ with $n_l \leq Ce^{Cl}$ and $n_l\rightarrow \infty$ satisfies 
 $$
  B_{\hat{g}_i}(\bar{x}_i, l)\subset \cup_{j=0}^{n_l} B_{\hat{g}_i}(x^{i}_j, 1/9) \subset B_{\hat{g}_i}(\bar{x}_i, l+1) $$
for $l=1,2\cdots$,  where $\Psi^i_j(0)=x^{i}_j$. 

Let $l(k_0)$ be the least $l\in \mathbb{N}$ with $m(k_0) \leq n_{l}$.  Then  $d_{\hat{g}_i}(\bar{x}_i, \Psi^i_{k_0}(y_{k_0})) \leq l(k_0)+1$ and  there exists  a $k^{\prime i}_0\leq n_{l(k_0)+1}\leq Ce^{Cl(k_0)+C}$ with  $\Psi^i_{k_0}(y_{k_0})\in B_{\hat{g}_i}(x^{i}_{k^{\prime i}_0}, 1/9)$.  After passing to a subsequence of $i$, we may assume all $k^{\prime i}_0\triangleq k^{\prime}_0$ are equal. 

 Hence there is a point $y^{i\prime}_{k_0} \in \Omega^{\prime}_{k^{\prime}_{0}}$ with $|y_{k^{i\prime}_0}|<1/8$ and  $\Psi^i_{k_0}(y_{k_0})=\Psi^i_{k^{\prime}_0}(y_{k^{i\prime}_0})$. 

Note that  $y_{k_0}\in \Omega^{\prime}_{m(k_0)}=\{|z|<1/2(1-\delta_0)\}$, there exist  a small $\xi_{k_0}>0$   with $\{z: |z-y_{k
_0}|<\xi_{k_0}\}\subset \Omega^{\prime}_{m(k_0)}$. Since $\Psi^i_{k_0}$ and $\Psi^i_{{k}^{\prime}_0}$ are local covering maps,  there exists  a map $L^{i}_{k_0}: \{z: |z-y_{k
_0}|<\xi_{k_0}\} \rightarrow \Omega^{\prime}_{m(k^{\prime}_0)}\triangleq \{|w|<1/2 (1-\delta_0)\}$ such that 
\begin{equation} \label{4.23}
\Psi^i_{k_0}=\Psi^i_{{k}^{\prime}_0} \circ L^{i}_{k_0}, \ \ \ \ L^i_{k_0}(y_{k_0})=y_{k^{i\prime}_0}.
\end{equation}
From (\ref{4.23}), we know  $L^{i}_{k_0}$ preserves $\hat{g}_i$, $\bar{g}_i$, $X_{i}$, $u_i$, $\omega_i$.   For example,  (\ref{4.23})  implies $$
(\Psi^i_{k_0})^{\ast} \hat{g}_{i}=  (L^{i}_{k_0})^{\ast}(\Psi^i_{k^{\prime}_0})^{\ast} \hat{g}_{i},$$
i.e.  $L^{i}_{k_0}$ is a  $\hat{g}_{i}$-isometry. Consequently,  
\begin{equation} \label{4.24}
\{|w-y_{k^{i\prime}_0}|<1/2\xi_{k_0}\}\subset L^i_{k_0}(\{z: |z-y_{k
_0}|<\xi_{k_0}\})\subset \{|w|<1/7\},
\end{equation}
because  $\xi_{k_0}$ is small and $|y_{k^{i\prime}_0}|<1/8$. Let $w^{\beta}= (L^{i}_{k_0}(z))^{\beta}$, then 

$$
\hat{g}_i(\partial_{z^{\beta}}, \partial_{z^{\gamma}})=\hat{g}_i(\partial_{w^{\xi}}, \partial_{w^{\eta}})\frac{\partial w^{\xi}}{\partial z^{\beta}}\frac{\partial w^{\eta}}{\partial z^{\gamma}}$$
and  
\begin{equation} \label{4.25}
\frac{\partial^2 w^{\beta}}{\partial z^{\gamma} \partial z^{\delta}}=\hat{\Gamma}^{\epsilon}_{\gamma\delta}(z)\frac{\partial w^{\beta}}{\partial z^{\epsilon}}-\hat{\Gamma}^{\beta}_{\xi\eta}(w) \frac{\partial w^{\xi}}{\partial z^{\beta}}\frac{\partial w^{\eta}}{\partial z^{\gamma}}.\end{equation}
(\ref{4.25}) follows from the standard computations of Christoffel symbols.  Since $\hat{g}_i\in W^{1,p_0}$,  we have $L^{i}_{k_0}\in W^{2,p_0}$ from (\ref{4.25}). By Sobolev imbedding theorem, $L^{i}_{k_0}\in C^{1,\alpha_0}$.
 
For any $0<\alpha<\alpha_0$, there is a subsequence of $i$ (assume we have already chosen the subsequence) such that $L^{i}_{k_0}\rightarrow L_{k_0}: \{|z-y_{k_0}|<\xi_{k_0}\}\rightarrow \{|w| <1/2(1-\delta_0)\}$ in $C^{1,\alpha}$ topology.  Denote $y^{\prime}_{k_0}=\lim_{i} y^{i\prime}_{k_0}\in \Omega^{\prime}_{m(k_0)}$, then $L_{k_0}(y_{k_0})=y^{\prime}_{k_0}$.  Since $\hat{g}_{i}, \bar{g}_{i}, X_{i}, u_{i}$ converge to $\hat{g}_{\infty}, \bar{g}_{\infty}, X_{\infty}, u_{\infty}$ in $C^{\alpha}$ topology, $L_{k_0}$ preserves $\hat{g}_{\infty}$, $\bar{g}_{\infty}$, $X_{\infty}$, $u_{\infty}$. Furthermore, this implies $L_{k_0}$ preserves $\omega_{\infty}$. 

  Consequently, we have 
 $$
 f(y_{k_0})=f(y^{\prime}_{k_0}).$$
 Hence 
\begin{equation} \label{4.26}
\begin{split}
& f(w)\leq \bar{l},  \ \ \ \text{on}  \ \Omega^{\prime}_{k_0^{\prime}}=\{|w|<1/2(1-\delta_0)\}\\
& f(y^{\prime}_{k_0}) \geq (1-k^{-1}_0) \bar{l}.
\end{split}
\end{equation}
Note that $|y^{\prime}_{k_0}|<1/7$.  Since for any $k\geq 1$,  the  $C^{k}$ norms (in $w$-coordinates) of $ (\hat{g}_{\infty}, \bar{g}_{\infty}, X_{\infty}, u_{\infty}, \omega_{\infty})$ on $\Omega^{\prime}_{m(k_0)}$ are  uniformly bounded,  there is a subsequence of $k_0\rightarrow \infty$ such that 
$ (\hat{g}_{\infty}, \bar{g}_{\infty}, X_{\infty}, u_{\infty}, \omega_{\infty})$ will converge to some limit (denoted by the same notation) in $C^{\infty}_{loc}$ topology. Denote $\lim\limits_{k_0\rightarrow \infty} y^{\prime}_{k_0}=\bar{y}$.   Clearly $|\bar{y}|\leq 1/7$, and 
\begin{equation} \label{4.27}
\begin{split}
& f(w)\leq \bar{f},  \ \ \text{on}  \ \{|w|<1/2(1-\delta_0)\},\\
& f(\bar{y}) =\bar{f}.
\end{split}
\end{equation}
In other words, $f$ achieves its maximum at an interior point $\bar{y}$. On the other hand,   the following Bochner type equation   holds  (cf. \cite{Chen1} (5.22),    \cite{A1}  (2.24)):   \begin{equation} \label{4.28}
\begin{split}
 \frac{1}{2}\triangle_{\hat{g}_{\infty}} f
 =  & 4|\nabla \log u_{\infty}|^4+|\omega_{\infty}|^2|\nabla \log u_{\infty}|^2 \\& +|2\nabla_{ij}\log u_{\infty} +(u_{\infty})^{-4}(\omega_{\infty})_i(\omega_{\infty})_j|^2 \\ & +(u_{\infty})^{-4}|\nabla_{i}(\omega_{\infty})_j-2(\omega_{\infty })_i(\log u_{\infty})_j-2(\omega_{\infty})_j(\log u_{\infty})_i|^2\\& +5|(u_{\infty})^{-2}\omega_{\infty}\wedge \nabla \log u_{\infty}|^2
\\& \  \geq 0.
 \end{split}
 \end{equation}
The strong maximum principle  implies that  $u_{\infty}\equiv 1$ and $\omega_{\infty}=0$, which is a contradiction with $f(\bar{y})=\bar{f}\neq 0$.   The claim is proved. Let $\delta=\delta_0\rightarrow 0$,  we still have $f\equiv 0$. This completes the proof. 
\end{pf}

\begin{rem}\label{r4.1} Let $\Psi^{k}_i:\Omega_i \rightarrow \bar{M}_k$ be the above immersion maps with $\cup^{n_k}_{i=1} \Psi^{k}_i(\Omega_i)\supset B_{\hat{g}_k}(\bar{x}_k, k)$.  Let  
$(G^1_{k}, G^{0}_k)$ be the groupoid 
which is equivalent to $\cup^{n_k}_{i=1} \Psi^{k}_i(\Omega_i)\subset \bar{M}_k$: 
\begin{equation} \label{4.29}
\begin{split} 
G^{0}_{k} & =\amalg^{n_k}_{i=0} \Omega_i,\\  G^{1}_{k} & =\amalg^{n_k}_{i,j=0} \{ (p_i,p_j): p_i\in \Omega_i, p_j\in \Omega_j, \Psi^{k}_i(p_i)=\Psi^{k}_j(p_j)\},
\end{split}
\end{equation} where the injection $e: G^{0}_k\rightarrow G^{1}_k$ is given by $e(p_i)=(p_i,p_i)$ for $p_i\in \Omega_i$; the source and  range maps  $s, r: G^1_k\rightarrow G^{0}_k$ are   $s((p_i,p_j))=p_j, r((p_i,p_j))=p_i$, the partially defined multiplication is the obvious one: $(p_i,p_j) (p_j,p_k)=(p_i,p_k)$.  Equip $(G^1_{k}, G^{0}_k)$ with the induced Lorentzian, Riemannian metrics $\bar{g}_k$, $\hat{g}_k$ and timelike vector field $X_k$ from $\bar{M}_k$, then the “arrows” in  $G^1_{k}$ become isometries preserving $\bar{g}_k$, $\hat{g}_k$, and $X_k$.    
The above proof based on (\ref{4.25}) essentially argues that the groupoids $(G^1_{k}, G^{0}_k)$ (around  some marked point $I^{i}_{k_0}$) have  a convergent subsequence in certain topology.  See \cite{H} \cite{Lott}  for the introduction of groupoids.    \end{rem}

\section{A  smoothing argument}
Analogous to the  weak $C^{\alpha}$ or $W^{1,p}$-harmonic radius for Lorentzian manifolds in Section 2,  the corresponding definition for Riemannian manifolds can be done by simply ignoring  the timelike vector fields.  Let us begin with some  terminologies, cf. \cite{PWY}, \cite{Pe}. 

Let $(M^n, g)$ be a Riemannian manifold, and $\Omega\subset M^n$ be an open subset. Fix $0<\alpha<1$, $r>0$, $k\in \mathbb{N}$, the weak harmonic $C^{k,\alpha}$-norm of $\Omega\subset M^{n}$ on scale $r$, denoted by $||\Omega||^{w,h}_{C^{k, \alpha}, r}$, is defined to be the infimum of the set of positive numbers  $Q$ such that there are a family of  local diffeomorphisms $\psi_{\tau}$ defined on a neighborhood of $B(0,r)\subset \mathbb{R}^n$:
\begin{equation} \label{5.1}
	\psi_{\tau}: B(0,r)(\subset \mathbb{R}^n)\rightarrow U_{\tau}(=\psi_{\tau}(B(0,r))\subset M^n), \ \ \   \tau\in \mathcal{I}, 
		\end{equation}
	satisfying 
	
	 $ 1) \ \  \triangle_{g_{\tau}} z^i=0,  \ \ \ \text{where} \ \ g_{\tau}\triangleq \psi^{\ast}_{\tau} g$;
	
	$ 2) \ \ \  e^{-2Q}\delta_{ij}\leq (g_{\tau})_{ij}\leq e^{2Q} \delta_{ij}$;
	
	3) \ for any $P\in \Omega$, there exists $\tau\in \mathcal{I}$, such that 
	$
	B(P,\frac{r}{10}e^{-Q})\subset U_{\tau}
	$;
	
	$ 4)   \ r^{k+\alpha}[(g_{\tau})_{ij}]_{C^{k,\alpha}(B(0,r))}\leq Q$.

	We remark that the definition of  $||\Omega||^{w,h}_{C^{k,\alpha},r}$ requires that the distance of  $\Omega$ to the ``boundary'' of $M^n$ is not less than $e^{-Q} r$.

If the map $\psi_{\tau}$ is required to be a diffeomorphism for all $\tau\in \mathcal{I}$, the norm is called harmonic $C^{k,\alpha}$-norm  (on scale $r$), denoted by $||\Omega||^{h}_{C^{k,\alpha},r}$.   Clearly,  one can define  (weak) harmonic $W^{k,p}$-norm by replacing  the $C^{k,\alpha}$-norm in 4) by

	$ 4)^{\prime}\ \ \ r^{|l|-\frac{n}{p}} ||\partial^l (g_{\tau})_{ij}||_{L^p(B(0,r))}\leq Q, \ \ \text{for} \ |l|\leq k.$ 
	
 Under the premise of the contradiction argument, the weak harmonic $C^{\alpha}$-norm of a big region  (in $(\bar{M}_j,\hat{g}_j)$) is bounded($\preceq \epsilon_0$) for scales $\leq 1/2$.  This is not enough for our later applications.   The purpose of this section is to prove that \textbf{the  weak harmonic $C^{\alpha}$-norm of the sequence   is arbitrarily small when $j>>1$}.  
 
 We have to take advantage of the behavior of the limit $(\bar{M}_{\infty}, \hat{g}_{\infty})$. 

For each $k$, $(\Omega_k, \hat{g}_{\infty})$ is flat,  $\Omega_k=\{|z|<1/2\}$, and the coordinates $\{z^{\beta}\}$ are harmonic for all $\hat{g}_j$ and $\hat{g}_{\infty}$. Let  $\{|z^{\prime}|<49/100\}$ be a rectangular coordinate system of $\hat{g}_{\infty}$ on $\Omega_k$ with the same origin as $\{z^{\beta}\}$ and 
\begin{equation} \label{5.2}
\{|z|\leq \frac{46}{100}\}\subset\{|z^{\prime}|\leq \frac{47}{100}\}\subset\{|z|\leq \frac{48}{100}\} \subset \{|z^{\prime}|<\frac{49}{100}\}. \end{equation}
This can be done if $\epsilon_0$ is suitably small in Section 2. 

Then  $z^{\prime\beta}(z)$ are  $C^{\infty}$ functions of $z=\{z^{\gamma}\}$. For each $j$, we compute 
$\triangle_{\hat{g}_j} z^{\prime\beta}$ in $z$-coordinates:

\begin{equation} \label{5.3}
\triangle_{\hat{g}_j} z^{\prime\beta}=(\hat{g}_j)^{\beta\gamma} \frac{\partial^2 z^{\prime\beta}}{\partial z^{\beta}\partial z^{\gamma}} = ((\hat{g}_j)^{\beta\gamma}-(\hat{g}_{\infty})^{\beta\gamma})\frac{\partial^2 z^{\prime\beta}}{\partial z^{\beta}\partial z^{\gamma}}.
\end{equation}
Since $\hat{g}_j\rightarrow \hat{g}_{\infty} $ in $C^{\alpha}$ topology (in $z$-coordinates), we know \begin{equation} \label{5.4}
||\triangle_{\hat{g}_j} z^{\prime\beta}||_{C^{\alpha}(|z| \leq \frac{48}{100})}\rightarrow 0, \ \ \ \text{as} \ \ j \rightarrow \infty.
\end{equation}

We solve the  Dirichlet problem:
\begin{equation} \label{5.5}
\begin{split}
 \triangle_{\hat{g}_j} w^{\prime \beta} & =-\triangle_{\hat{g}_j} z^{\prime \beta}    \ \ \text{in} \ \{|z|<\frac{48}{100}\}\\
  w^{\prime \beta} & =0  \ \ \text{on}  \ \{|z|=\frac{48}{100}\}.
\end{split}
\end{equation}

Let $w^{\beta}=z^{\prime\beta}+w^{\prime \beta}-w^{\prime \beta}(0)$, then $w^{\beta}$ is $\hat{g}_{j}$-harmonic, and 
\begin{equation} \label{5.6}
|| w^{\prime\beta} ||_{C^{2,\alpha}(|z|\leq \frac{48}{100})}\leq C||\triangle_{\hat{g}_j} z^{\prime\beta}||_{C^{\alpha}(|z| \leq \frac{48}{100})}\rightarrow 0\end{equation}
by Schauder estimates in $z$-coordinates. Since the coordinate  transformation $\{z^{\beta}\}\rightarrow \{z^{\prime \beta}\}$ is $C^{\infty}$-smooth, we have 
\begin{equation} \label{5.7}
|| w^{\prime\beta} ||_{C^{2,\alpha}(|z^{\prime}|\leq \frac{47}{100})}\rightarrow 0\end{equation}
where the $C^{2,\alpha}$ norm is in $z^{\prime}$-coordinates. So $\{w^{\beta}\}$ is a coordinate system for large $j$. 
By same reason,  $\hat{g}_j(\frac{\partial}{\partial z^{\prime \gamma}}, \frac{\partial}{\partial z^{\prime \delta}})\rightarrow \hat{g}_{\infty}(\frac{\partial}{\partial z^{\prime \gamma}}, \frac{\partial}{\partial z^{\prime \delta}})=\delta_{\gamma\delta} $ in $C^{\alpha}$ topology (in $z^{\prime}$-coordinates). Combining with (\ref{5.7}),  
we have 
\begin{equation} \label{5.8}
||\hat{g}_j(\frac{\partial}{\partial w^{ \gamma}}, \frac{\partial}{\partial w^{\delta}})-\delta_{\gamma\delta} ||_{C^{\alpha}(|w|\leq \frac{46}{100})} \rightarrow 0
\end{equation}
as $j\rightarrow \infty$. 

 Recall that 
for any $0<l<1/4 m_j \rightarrow \infty$, 
$\{B_{\hat{g}_j}(x^{j}_k, \frac{1}{9})\}_{k=0}^{n_l}$ is a covering of $B_{\hat{g}_j}(\bar{x}_{j}, l)$. If we choose $r=\frac{1}{3}$, for any point $P\in B_{\hat{g}_j}(\bar{x}_{j}, l)$, there is a $k\leq n_l$ with  $P\in B_{\hat{g}_j}(x^{j}_k, \frac{1}{9})$, hence $B_{\hat{g}_j}(P, \frac{1}{10} r)\subset B_{\hat{g}_j}(x^{j}_k, \frac{1}{9}+\frac{1}{30}) \subset \Psi^{j}_{k}(|w|<\frac{1}{3})$. In particular, this implies 
\begin{equation} \label{5.9}
||B_{\hat{g}_j}(\bar{x}_{j}, l) ||_{C^{\alpha}, r}^{w,h}\rightarrow 0
\end{equation}
as $j\rightarrow \infty$, where $r=\frac{1}{3}$. 

Now we try to smooth the Riemannian metric $\hat{g}_j$ a little, so that the resulting metric  is close to $\hat{g}_j$ and has very small curvatures. To this end, we shall  use the construction  in \cite{PWY}. The following  theorem is essentially a generalization of Theorem 1.1 in \cite{PWY}.  

\begin{thm} \label{t5.1}  For any $\epsilon, \delta, \alpha \in (0,1)$, $0<\alpha^{\prime}<\alpha$,  $n,k\in \mathbb{N}$, there is a small constant $Q(\delta,\alpha, \epsilon, n, k)>0$ such that if $(M^n, g)$ is  a Riemannian manifold, $\Omega\subset M^n$ is an open subset, satisfying
	\begin{equation} \label{5.10}
		||(\Omega,g)||^{w,h}_{C^{\alpha}, r_0}\leq Q,
				\end{equation} 			 then there exist a new Riemannian metric $g_{\epsilon}$ and a constant $\tilde{Q}\in (0,\delta)$  satisfying 
\begin{equation} \label{5.11}
e^{-2\epsilon} g\leq g_{\epsilon} \leq e^{2\epsilon} g
\end{equation}
\begin{equation} \label{5.12}
	\begin{split}
	\sum_{l=0}^{k}r^{2+l}_0| \nabla_{g_{\epsilon}}^{l} Rm(g_{\epsilon})|_{\Omega} \leq \tilde{Q} 
		\end{split}
	\end{equation}	
	and 
\begin{equation} \label{5.13}
\sup_{\tau\in I}[g_{\epsilon}]_0+ [g_{\epsilon}]_{1, \alpha^{\prime}}  r^{1+\alpha^{\prime}}_0\leq \tilde{Q},  
\end{equation}
  where the supremum is taken over  the family of $g$-harmonic coordinates $\psi_{\tau}: B(0,r_0)\rightarrow U_{\tau}=\psi_{\tau}(B(0,r_0))\subset M^n$, $\tau\in I$, $\cup_{\tau\in I} U_{\tau} \supset \Omega$,   in the definition of $||(\Omega, g)||^{w,h}_{C^{\alpha}, r_0}$.
  \end{thm}

We remark that   for any point $P'\in \Omega$, the metric  $g_{\epsilon}(P')$ in Theorem \ref{t5.1} depends only on the metric $g$ on the ball $B(P', 2r_0)$.   The dependence of $\tilde{Q}$ on $Q$, i.e.,  $\lim\limits_{Q\rightarrow 0} \tilde{Q}=0$,  is important  for later applications. To this end, we need the  concrete construction in \cite{PWY}.

 Let $\psi_{\tau}:B(0,r_0)\rightarrow U_{\tau}$ be a local diffeomorphism in the definition of weak harmonic $C^{\alpha}$-norm on scale $r_0$. Let $\Omega'=B_{g}(0, (1-10^{-2}) r_0)\subset B(0,r_0)$ be a geodesic ball of radius $(1-10^{-2}) r_0$ in the pulled back metric $g$, $\Omega'_1=\Omega'- \cup _{q\in \partial \Omega'} \bar{B}_{g}(q, \frac{r_0}{10})$. For $i_0= r_0/ 10$, and $x\in \Omega'_1$,   solve the Dirichlet problem:
	\begin{equation} \label{5.14}
		\begin{split}
		 \triangle_{g} h_{x} & =-1 \ \  \text{in} \ B_{g}(x, i_0)\\
		 h_{x} & = 0 \ \ \text{on} \ \partial B_{g}(x, i_0).
		\end{split}
		\end{equation}
	
	Let	
	$\tilde{\beta}:\mathbb{R} \rightarrow \mathbb{R}$
	be a smooth nonnegative nondecreasing function with  $\tilde{\beta}(t)\equiv 0$ for $t\leq \frac{1}{2}$ and $\tilde{\beta}(t)\equiv 1$ for $t \geq 2$.  Define 
	\begin{equation} \label{5.15}
		f(p,q)=\int_{\Omega'_1} \tilde{\beta}(\frac{2n}{i^2_0} h_{x}(p))\tilde{\beta}(\frac{2n}{i^2_0} h_{x}(q)) dv_{g}(x),
		\end{equation}
	for $p,q\in \Omega'$. 
	When $g$ is  Euclidean, $\bar{h}_{x}(\cdot)=\frac{1}{2n}(i^2_0-d^2_{E}(x,\cdot))$ is the  solution to (\ref{5.14}). In this case, the corresponding $\bar{f}(p,q)=i_0^n\tilde{f}(i_0^{-1} d_{E}(p,q))$, where $\tilde{f}$ is a smooth  function depending only on $\tilde{\beta}$ and $n$.  Let  $\tilde{F}: \Omega'_1 \rightarrow L^{2}(\Omega^{\prime}, g_{E})$ be  defined by $p\rightarrow i_0^{-\frac{3 n}{2}+1} \bar{f}(p,\cdot)$. By direct computations, one can prove 
	$$
	||d_{v_p}\tilde{F}||^2_{L^2(\Omega_1)}=[\frac{vol(S^{n-1})}{n} \int_{0}^2 r^{n-1} \tilde{f'}^2(r)dr] |v_p|^2= B^{-2}_n |v_p|^2
	$$
where $B_n=[\frac{vol(S^{n-1})}{n} \int_{0}^2 r^{n-1} \tilde{f'}^2(r)dr]^{-\frac{1}{2}}$. Let $\beta(t)=B_n \tilde{\beta}(t)$, and 
\begin{equation} \label{5.16}
	\begin{split} 
	f_{p}(q)&=\int_{\Omega'_1} {\beta}(\frac{2n}{i^2_0} h_{x}(p)){\beta}(\frac{2n}{i^2_0} h_{x}(q)) dv_{g}(x)\\
	F(p)&=i_0^{-\frac{3n}{2}+1} f_p: \Omega^{\prime}_1\rightarrow L^2(\Omega^{\prime}, g).
	\end{split}
	\end{equation}
Then  $F$ is an isometric embedding when $g$ is Euclidean.

Direct computation gives ((4.11) and (4.12) in \cite{PWY})

\begin{equation} \label{5.17}
	\partial_{z^i} F(p)=2n i_0^{-\frac{3}{2}n} \int_{\Omega_1}\beta'(\frac{2n}{i_0^2} h_x(p)) i_0^{-1} \partial_{z^i} h_{x}(p)\beta(\frac{2n}{i_0^2} h_x(q))dv_g(x)
\end{equation}
\begin{equation} \label{5.18}
	\begin{split}
		\partial^2_{z^i z^j} F(p) &=4n^2 i_0^{-\frac{3}{2}n-1} \\ & \times \int_{\Omega_1}\beta''(\frac{2n}{i_0^2} h_x(p)) i_0^{-1} \partial_{z^i} h_{x}(p) i_0^{-1} \partial_{z^j} h_{x}(p) \beta(\frac{2n}{i_0^2} h_x(q))dv_g(x)\\
		&	+2n i_0^{-\frac{3}{2}n-1} \int_{\beta_1}\beta'(\frac{2n}{i_0^2} h_x(p))  \partial^2_{z^i z^j} h_{x}(p)\beta(\frac{2n}{i_0^2} h_x(q))dv_g(x).
	\end{split}
\end{equation}

Let  ${G}=i_0^{-2} g$, then  $B_{g}(x,i_0)$ becomes $B_{G}(x,1)$.   Let   $w^i=i^{-1}_0 z^i$, $G(\frac{\partial}{\partial w^i},\frac{\partial}{\partial w^j})\triangleq G_{ij}=g_{ij}$, $H_{x}(w)=i_0^{-2}h_{x}(i_0w)$. 	In harmonic coordinates $\{w^{i}\}$, the
equation (\ref{5.14}) becomes  
\begin{equation} \label{5.19}
\triangle_{G} H_x=G^{ij} \partial^2_{w^i w^j} H_x=-1.
\end{equation}

Multiplying (\ref{5.19}) by $H_x$,  integrating by parts and using Poincar$\acute{e}$ inequality, we get

\begin{equation} \label{5.20}
		\int_{B_{G}(x,1)}  |\nabla H_x|^2+H^2_x dw \leq C.
\end{equation}

Since $||G_{ij}-\delta_{ij}||_{C^{\alpha}(|w|<1)}\leq Q$, combining (\ref{5.19}) and (\ref{5.20}),   we have the uniform interior $C^{2,\alpha}$ estimate for  $H_{x}$ i.e. for any $0<R<1$, there is a constant $C(Q, R,\alpha)$  such that 
\begin{equation} \label{5.21}
	\begin{split}
||H_x||_{C^{2,\alpha}(|w|\leq R)}\leq C(Q, R, \alpha).
\end{split}
	\end{equation}

Denote $\bar{H}_{x}(w)=\frac{1}{2n}(1-|w|^2)$.   When $Q\rightarrow 0$, as observed  in \cite{PWY},   $H_x$ will weakly converge to $\bar{H}_{x}$ in $W^{1,2}$,  and 
\begin{equation} \label{5.22}
||H_x-\bar{H}_x||_{C^{2, \alpha^{\prime}}(|w|<R)}\rightarrow 0,
\end{equation} for $0<R<1$,  $0<\alpha^{\prime} <\alpha$. 
Note that 
\begin{equation} \label{5.23}
\begin{split}(F^{\ast} g_{L^2})_{ij}(p)&=4n^2\int_{\Omega'\times \Omega'} \beta^{'}(2nH_{x}(p))\beta^{'}(2nH_{y}(p))\partial_{w^i} H_{x}(p)\partial_{w^j}H_{y}(p)\\ 
		&\ \ \ \times B(x,y) dv_{G}(x) dv_{G}(y)
		\end{split}
\end{equation}
 where $$B(x,y)=\int_{\Omega'}\beta(2n H_x(q))\beta(2n H_y(q))dv_{G}(q).$$ 
 
From (\ref{5.21}) and (\ref{5.23}), we know the $C^{1,\alpha}$ norm of  $(F^{\ast} g_{L^2})_{ij}$ is uniformly bounded. Moreover, (\ref{5.22})  implies
  \begin{equation} \label{5.24}
  ||(F^{\ast} g_{L^2})_{ij}-\delta_{ij}||_{C^{1, \alpha^{\prime}}(|w|<R)}\rightarrow 0, 
  \end{equation} for $0<\alpha^{\prime}<\alpha$,  as $Q\rightarrow 0$. 

On the other hand,  the  Gauss equation gives 
\begin{equation} \label{5.25}
	Rm( F^{\ast} g_{L^2})_{ijkl}=\langle \nabla_i\nabla_k F, \nabla_j\nabla_lF \rangle_{L^2}-\langle \nabla_i\nabla_l F, \nabla_j\nabla_kF \rangle_{L^2}
\end{equation}
where $$\nabla_i\nabla_j F=\partial^2_{z^iz^j}F-
\langle\partial^2_{z^iz^j}F,\partial_{z^l} F\rangle_{L^2} (F^{\ast}g_{L^2})^{kl} \partial_k F.$$

Similar to (\ref{5.23}), the curvature $Rm( F^{\ast} g_{L^2})_{ijkl}$ has the expression of the form: 
\begin{equation}
	\begin{split} \label{5.26}
& [\int_{\Omega'\times\Omega'}  (\beta^{''} \partial H \ast \partial H+\beta'\partial^2 H)(x)(\beta^{''} \partial H \ast \partial H+\beta'\partial^2 H)(y)\\ & \times B(x,y)dv_{G}(x)dv_{G}(y)+ (\int_{\Omega'\times\Omega'}  (\beta^{''} \partial H \ast \partial H+\beta'\partial^2 H)(x)(\beta^{'} \partial H)(y)\\ & \times B(x,y)dv_{G}(x)dv_{G}(y))^2 \ast (F^{\ast} g_{L^2})^{-1}] i^{-2}_0 \\
&\triangleq E(\beta, \beta',\beta'', H, \partial H, \partial^2H, G) i^{-2}_0.
	\end{split}
	\end{equation}
 (\ref{5.21}) shows that the $C^{\alpha}$ norm of $E(\beta, \beta',\beta'', H, \partial H, \partial^2H, G)$ is uniformly bounded.  Moreover, as $Q\rightarrow 0$, (\ref{5.22}) implies that  $E(\beta, \beta',\beta'', H, \partial H, \partial^2H, G)$ converges in $C^{\alpha^{\prime}}$-topology (for $0<\alpha^{\prime}<\alpha$) to 
 $$
 E(\beta,\beta',\beta'', \bar{H}, \partial{\bar{H}}, \partial^2{\bar{H}}, \delta)\equiv 0,
 $$
 which is the curvature of the  Euclidean metric by our construction.  
 
 It is crucial  to observe that the above construction is equi-variant, i.e.  $F^{\ast}g_{L^2}$ is also invariant under $g$-isometries. Hence  $F^{\ast}g_{L^2}$ can descend to the manifold $U_{\tau}$,  and coincides on any  intersections $U_{\tau}\cap U_{\tau^{\prime}}$.  In other words, we obtain a globally defined metric on $\Omega$.  
 
 \begin{pf} of Theorem \ref{t5.1}. 
 
   For $k=0$, the above construction gives the desired metric $(g)_{\epsilon}\triangleq g_1$.  (\ref{5.11}) and   (\ref{5.13}) follow from (\ref{5.24}).   Since  $|Rm(g_1)|r^{-2}_0\rightarrow 0$ as $Q\rightarrow 0$,  by \cite{JK}, its weak $C^{1,\alpha}$-harmonic radius is not less than  $r_0$, more precisely, we have  $||(\Omega,g_1)||^{w,h}_{C^{1,\alpha}, r_0}\rightarrow 0$, as $Q\rightarrow 0$.  Repeating the above construction to smooth $g_1$ again, one obtain a new metric $(g_{1}){\epsilon}\triangleq g_2$.  Since $||G^{ij}-\delta_{ij}||_{C^{1,\alpha}}<\tilde{Q}$  in (\ref{5.19}),  (\ref{5.24}) (in $g_1$-harmonic coordinates) is now replaced by:  
 \begin{equation}  \label{5.27}
  [(F^{\ast} g_{L^2})_{ij}-\delta_{ij}]_0+[(F^{\ast} g_{L^2})_{ij}-\delta_{ij}]_{2, \alpha^{\prime}} r^{2+\alpha^{\prime}}_0\rightarrow 0.
 \end{equation}
 Taking one more derivative on (\ref{5.25}) and (\ref{5.26}),   we have 
  \begin{equation} \label{5.28} 
  |Rm((g_1)_{\epsilon})| r^{-2}_0+|\nabla Rm((g_1)_{\epsilon})|r^{-3}_0\rightarrow 0.
 \end{equation} 
  Since (\ref{5.13}) holds for $g_1$, by Schauder estimates, $g_1$-harmonic coordinates are $C^{2,\alpha^{\prime}}$ bounded in $g$-harmonic coordinates. Let $Q\rightarrow 0$,  using (\ref{5.27}), we know  (\ref{5.13}) holds for $g_2$. Repeating the above argument  many times, we obtain a sequence $g_1, g_2,\cdots, g_{k+1}$ and $g_{k+1}$ satisfies (\ref{5.13}) and  (\ref{5.12}).   
  \end{pf}\\
 As a consequence,    one can prove:
\begin{prop} \label{p5.2}
For large $j$, one can  construct  a smooth non-degenerate map $\Phi_j:\{ |w|<a_j\}\rightarrow \bar{M}_j$ with $\Phi(0)=\bar{x}_j$, $a_j\rightarrow \infty$ as $j\rightarrow \infty$,   and let $\hat{g}_{\alpha\beta}\triangleq (\Phi_j^{\ast} \hat{g}_j)(\frac{\partial }{\partial w^{\alpha}}, \frac{\partial }{\partial w^{\beta}})$,  we have

\begin{equation} \label{5.29}
	\frac{1}{2} \delta_{\alpha\beta} \leq  \hat{g}_{\alpha\beta}\leq 2 \delta_{\alpha\beta}, \ \ \ \text {on} \ \  \{|w|<a_j\}.
\end{equation}

\end{prop}
\begin{pf}
By (\ref{5.9}),  for any $l>0$, the weak harmonic $C^{\alpha}$-norm $||B_{\hat{g}_j}(\bar{x}_j, l)\subset \bar{M}_j||^{w,h}_{C^{\alpha}_{, r}}$  with $r=1/3$ converges to zero, as $j\rightarrow \infty$. By Theorem \ref{t5.1}, there is a metric $g_{\epsilon}$ on $B_{\hat{g}_j}(\bar{x}_j, l)$ such that 
$$
e^{-2\epsilon} \hat{g} \leq g_{\epsilon} \leq e^{2\epsilon} \hat{g}
$$
and $|Rm(g_{\epsilon})|\rightarrow 0$ as $j\rightarrow \infty$. Let $\exp_j: T_{\bar{x}_j}\bar{M}_j\rightarrow \bar{M}_j$ be the exponential map of the metric $g_{\epsilon}$, $I:\mathbb{R}^n \rightarrow T_{\bar{x}}\bar{M}_j$ be a linear isometry, $\Phi_j=\exp_j \circ I$.   (\ref{5.29}) follows from the standard Rauch comparison theorem  of Jacobi fields in Riemannian geometry.
\end{pf}

\section{Harmonic radius}

	Let $\Phi_j:\{ |w|<a_j\}\rightarrow \bar{M}_j$ be the map in Proposition \ref{p5.2}. We  use $\Phi_j$ to pull back  $\hat{g}_j$, $\bar{g}_j$, $X_j$ from $\bar{M}_j$ to $\{ |w|<a_j\}$  and  maintain   the same notations for them.  As in Remark \ref{r2.1}, one can define 	$W^{1,p}$-harmonic radius by  modifying Definition \ref{d2.1} in Section 2 to require the immersion map to be a diffeomorphism.

	Denote  $Q_j=\{ w: d_{\hat{g}_j}(w,0)<\frac{a_j}{2}\}\subset \{|w|<a_j\}$. 
The proof of Theorem \ref{t1.1} is reduced  to the following proposition. 
	\begin{prop} \label{p6.1} For any $\epsilon\in (0,1)$ and $p>4$, there is a constant $C$ depending only on $\epsilon$ and $p$ such that  on the manifold $(Q_j, \bar{g}_j, \hat{g}_j, X_j)$, we have 
		$$\sup_{Q_j} (\frac{a_j}{2}-d_{\hat{g}_j}(w, \partial Q_j) )(r^{h}(W^{1,p}, \epsilon)(w)) ^{-1}\leq C.$$
		\end{prop}

	\begin{pf}
		We argue by contradiction as in Section 2. Suppose the result is not true,  there exist  $\epsilon_1>0$,  $p_1>4$ and  a sequence of pointed Lorentzian manifolds  with timelike vector fields $(Q_j, \bar{g}_j, X_j, O_j)$, such that  the rescaling of $\bar{g}_j, \hat{g}_j, X_j$ (denote them with the same notations as before) satisfies 
\begin{equation} \label{6.1}
	\begin{split} |X_j(\bar{w}_j)|^2_{\hat{g}_j}& =1, \ \ \  r^{h}(W^{1,p_1}, \epsilon_1)(\bar{w}_j) =1,\\ 
	r^{h}(W^{1,p_1}, \epsilon_1)({x})& \geq \frac{1}{2}, \ \ \ \  \text{on} \ \   B_{\hat{g}_j}(\bar{w}_j, b_j),	\end{split}\end{equation}	
	and 	
	 \begin{equation}\begin{split} \label{6.2}
 \sup_{{B}_{\hat{g}_j}(\bar{w}_j, b_j)} |{R}ic(\bar{g}_j)|_{\hat{g}_j}  \rightarrow 0, \ \ \ \ 
 \sup_{{B}_{\hat{g}_j}(\bar{w}_j, b_j)} u_j^{-1}|\mathcal{L}_{X_j}\bar{g}_j|_{\hat{g}_j} \rightarrow 0,
  \end{split}
\end{equation}	
 where ${B}_{\hat{g}_j}(\bar{w}_j, b_j)\subset Q_j$, $b_j\rightarrow \infty$.

 Now,  the sequence of pointed Riemannnian manifolds  $({B}_{\hat{g}_j}(\bar{w}_j, b_j), \hat{g}_j, \bar{w}_j)$ has   a subsequence converging  to some complete Riemannnian manifold  $(M^{\infty}, \hat{g}^{\infty}, \bar{w}^{\infty})$ in $C^{\alpha}$ and weak $W^{1,p_1}$ topologies in Cheeger-Gromov sense, where $0<\alpha<\alpha_1=1-4/p_1$.  That is to say, there is a sequence of diffeomorphisms  $\phi_j$ from open subsets of $M^{\infty}$ to ${Q}_j$ with $\phi_j(\bar{w}_{\infty})=\bar{w}_j$ and  $\phi^{\ast}_j \hat{g}_j$ is   converging  to $\hat{g}^{\infty}$ in $C^{\alpha}$ and weak-$W^{1,p_1}$ topologies  on any  compact sets of $M^{\infty}$.  Note that  we  also need to investigate the convergence of  $\phi^{\ast}_j X_j$ and $\phi^{\ast}_j \bar{g}_j$.  For the sake of completeness,  we shall revisit   the  concrete construction of the maps  $\phi_j$, cf.  Chapter 10 in \cite{Pe}.

 First of all, $\{({B}_{\hat{g}_j}(\bar{w}_j, b_j), \hat{g}_j, \bar{w}_j)\}^{\infty}_{j=1}$ has a subsequence Gromov-Hausdorff converging to some complete pointed limit  space $(L, d_{L}, \bar{w})$. Recall that if $(X_i, d_i)$, $i\in I$ is a family of metric spaces, a metric $d$ on the disjoint union $\coprod_{i\in I} X_i$ is said to be compatible if $d|_{X_i}=d_{X_i}$. By selecting a further subsequence, one can find compatible metrics $d_{i,i+1}$ on $\bar{B}(\bar{w}_i, i)\coprod \bar{B}(\bar{w}_{i+1}, i)$, where $\bar{B}(\bar{w}_i, i)\subset {Q}_i$,  $\bar{B}(\bar{w}_{i+1}, i)\subset {Q}_{i+1}$ are closed balls of radius $i$, such that $d(\bar{w}_i, \bar{w}_{i+1})\leq 1/2^{i+1}$, and 
   \begin{equation} \label{6.3}
 d_{Hausdorff} (\bar{B}(\bar{w}_i, i),\bar{B}(\bar{w}_{i+1}, i) ) \leq 1/ 2^{i+1}.
 \end{equation}
 Extend $d_{i,i+1}$ to a compatible metric to $\bar{B}(\bar{w}_i, i)\coprod \bar{B}(\bar{w}_{i+1}, i+1)$ 
 by setting 
   \begin{equation} 
   \label{6.4}
    d_{i,i+1}(x_i,x_{i+1})=\inf_{p\in \bar{B}(\bar{w}_{i+1}, i)} d(x_i,p)+d(p,x_{i+1}),
 \end{equation}  
 for $x_i\in \bar{B}(\bar{w}_i, i)$, $x_{i+1}\in \bar{B}(\bar{w}_{i+1}, i+1)$.

 Define a compatible metric $d_{i,i+j}$ on  $\bar{B}(\bar{w}_i, i)\coprod \bar{B}(\bar{w}_{i+j}, i+j)$ 
 by setting 
  \begin{equation} \label{6.5}
 d_{i,i+j}(x_i,x_{i+j})=\inf_{x_{i+k}\in \bar{B}(\bar{w}_{i+k}, i+k)} d(x_i,x_{i+1})+\cdots+d(x_{i+j-1},x_{i+j}),
 \end{equation}
  for $x_i\in \bar{B}(\bar{w}_i, i)$, $x_{i+j}\in \bar{B}(\bar{w}_{i+j}, i+j)$.  
  
  This gives a compatible metric $d$ on $\coprod_{i}\bar{B}(\bar{w}_i, i)$. Adding the limit space $(L, d_{L}, \bar{w})$ to the ``boundary'' of the completion of $(\coprod_{i}\bar{B}(\bar{w}_i, i), d)$, we obtain a compatible metric (denoted by $d_{Y}$) on 
  $$
  Y=L \coprod \coprod_i \bar{B}(\bar{w}_i,i)
  $$ 
  satisfying $d_{Y}(\bar{w}_i, \bar{w}_{i+j})\leq 1/2^{i}$, $d_{Y}(\bar{w}_i, \bar{w}) \leq 1/ 2^{i}$ and 
  \begin{equation} \label{6.6}
  \begin{split}
  & d_{Hausdorff} (\bar{B}(\bar{w}_i, i), \bar{B}(\bar{w}_{i+j}, i)) \leq 1/2^{i}\\
   & d_{Hausdorff} (\bar{B}(\bar{w}_i, i), \bar{B}_{L}(\bar{w}, i)) \leq 1/2^{i}.  \end{split}
  \end{equation}
  
    As in Section 2,  there are   a  sequence 
of natural numbers $n_1\leq n_2\leq n_3\cdots$  with $n_l\leq Ce^{Cl}$ and a sequence of points $w_0=\bar{w}, w_1, w_2 \cdots $ in ${L}$ such that   $\{w_i\}_{i=0}^{n_l}$ is a set of maximal number of points in $B_{{L}}(\bar{w}, l)$ with mutual distances $\geq 1/20$. 

 For fixed $l$, let $\{x_{i}^{l}\}_{i=0}^{n_l} \subset B(\bar{w}_l,l) \subset \bar{M}_l$ satisfy 
 $$
 d_{Y}(x^{l}_{i},w_i) \leq 1/2^{l}, \ \ \ \text{for}\ \ i=0,1,\cdots n_l.
 $$
 Then $\{B(x^{l}_i, \frac{1}{40})\}_{i=0}^{n_l}$ is disjoint and $\{B(x^{l}_i, \frac{1}{19})\}_{i=0}^{n_l}$ is a covering of $B(\bar{w}_l, l) \subset {Q}_l$. 
    Since $r^{h}(W^{1,p_1}, \epsilon_1)(x^{l}_i) \geq 1/2$, for all $0\leq i \leq n_l$, there exist a family of diffeomorphisms  $\Psi^{l}_i: \{ |z|<1/2\} \rightarrow \bar{M}_l$  with $\Psi^{l}_i(0)=x^{l}_i$ and (\ref{2.1}) (\ref{2.2}) (\ref{2.3})  hold for $r=1/2$.

    Note that any bounded set of  $(Y, d_{Y})$ is paracompact. This fact  can be proven as follows. Let $\{y_k\}_{k=1}^{\infty}\subset Y$ be a bounded sequence, from which we will extract a Cauchy subsequence. 
   Let $y_k\in {Q}_{i_k}$ or $L$, one may regard $i_{k}=\infty$ for the latter case. The case when  $\{i_k\}_{k=1}^{\infty}$ has a bounded subsequence is trivial. So we may assume $i_k\rightarrow \infty$ as $k\rightarrow \infty$. Since $d_{Y}(\bar{w}_i,\bar{w}_j)\leq 1/2^{min\{i,j\}}$, there is an integer $D>0$ such that $d_{Y}(y_k,\bar{w}_i)\leq D$ for all $i, k$.   On ${Q}_{D}$, there is a bounded sequence $\{y_k^{\prime D}\}_{i_k\geq D}$ satisfying $d_{Y}(y_k, y_k^{\prime D})\leq 1/2^{D}$. Since $\{y_k^{\prime D}\}_{i_k\geq D}$ has a convergence subsequence in   $\bar{M}_{D}$, we know $\{y_k\}^{\infty}_{k=1}$ has a subsequence $\{y^{D}_k\}^{\infty}_{k=1}$ such that $d_{Y}(y^{D}_l, y^{D}_{l^{\prime}})\leq 3/2^{D}$. Inductively, one can construct subsequences 
   $$
   \{ y_k\}_{k=1}^{\infty}\supset  \{ y^{D}_k\}^{\infty}_{k=1} \supset  \{ y^{D+1}_k\}^{\infty}_{k=1}\supset \cdots  $$
   with 
     $d_{Y}(y^{D+j}_l, y_{l^{\prime}}^{ D+j})\leq 3/2^{D+j}$ for all $l,l^{\prime} >0$. The diagonal process gives the final Cauchy subsequence $\{y^{D+k}_k\}_{k=1}^{\infty}$.

     Note that for fixed $i$, the maps $\Psi^{k}_i:\{|z|\leq {7}/ {16}\}\rightarrow {Q}_k\subset Y$ are equi-continuous.  Since any bounded set of $Y$ is paracompact (this fact is crucial), we know $\{\Psi^{k}_i\}_{k=1}^{\infty}$ has  a $C^{0}$ convergent subsequence.  A diagonal process gives rise to a subsequence  $\{\Psi^{k_p}_i \}_{p=1}^{\infty} $ which converges for each $i$ in  $C^{0}$ to a  continuous limit map  $\Psi_i: \{|z|\leq 7/16\}\rightarrow L$ with $\Psi_i(0)=x_i$. 
     
     The limit  $\Psi_i$   still satisfies 
     
     $$ e^{c}|z-z^{\prime}|\geq d(\Psi_{i}(z), \Psi_{i}(z^{\prime})) \geq e^{-c}|z-z^{\prime}|$$
     for any $z,z^{\prime}\in \{|z|\leq 7/16\}$, where $c$ is a universal constant. Let $\Omega_i=\{|z|<7/16\}$,  then $\Psi_i: \Omega_i\rightarrow \Psi_i(\Omega_i)$ is a homeomorphism and   $\{\Psi_i (\Omega_i)\}_{i=1}^{\infty}$  forms  an open covering of $L$.   These  coordinate charts at least  make  $L$  to be a topological manifold. 
     
     On each ${Q}_{k}$,  the transition function $\Psi^{k}_{ij}= \Psi^{k}_i(\Psi^{k}_j)^{-1}$, from $\Omega_j$ to $\Omega_i$,   satisfies 
	\begin{equation} \label{6.7}
	\Psi^{k}_{li}\Psi^{k}_{ij}=\Psi^{k}_{lj}	\end{equation}		
	if $\Psi^{k}_i (\Omega_i) \cap 	\Psi^{k}_j (\Omega_j) \cap 	\Psi^{k}_l (\Omega_l) \neq \phi$. 
	
	 As in (\ref{4.25}) in Section 4,  $\Psi^{k}_{ij} \in W^{2,p_1}((\Psi^{k}_j)^{-1}\Psi^{k}_i (\Omega_i) \cap \Omega_j)$ and  it’s  $W^{2,p_1}$ norm  is uniformly bounded (independent of $k$).  Note that these transition functions $\Psi^{k}_{ij}$ may not be defined on same domains for different $k$. Suppose 	 $\Psi_i (\Omega_i) \cap \Psi_j (\Omega_j)\neq \phi$ for the limit,  then the  domains $(\Psi^{k}_j)^{-1}\Psi^{k}_i (\Omega_i) \cap \Omega_j$ will converge to some limit (denoted   by $\Omega^{\infty}_{ij}$), and $\Psi^{k}_{ij}$ has a subsequence which 	 converges to $(\Psi_j)^{-1} \Psi_i$ in $C^{1,\alpha}$ and weak $W^{2,p_1}$ topologies  on any compact subsets of $\Omega^{\infty}_{ij}$. 	 Using a diagonal process, one may assume the subsequence is chosen so that the transition functions are convergent for any index  pair $(i,j)$ with $\Psi_i (\Omega_i) \cap \Psi_j (\Omega_j)\neq \phi$.  The limit transition functions $(\Psi_j)^{-1} \Psi_i$ are  in  $W^{2,p_1}$, and  the cocycle condition (\ref{6.7}) also holds.

       There is a subsequence of   $\hat{g}_{k}\triangleq (\Psi^{k}_i)^{\ast}\hat{g}_{k} $, $\bar{g}_k \triangleq (\Psi^{k}_i)^{\ast}\bar{g}_{k}$ and $X_k\triangleq (\Psi^{k}_i)^{\ast}X_{k}$ which converges  to $\hat{g}_{\infty i}$, $\bar{g}_{\infty i}$ and $X_{\infty i}$ in $C^{\alpha}$ and weak $W^{1,p_1}$ topologies on $\{|z|< 7/16\} \triangleq \Omega_i$.  Using a diagonal process, one may assume the subsequence is chosen so that it converges on any $\Omega_i$ for all $0\leq i <\infty$.

Since $(\Psi^k_i)^{\ast} \hat{g}_k= ((\Psi^k_j)^{-1}\circ \Psi^k_i)^{\ast} (\Psi^k_j)^{\ast} \hat{g}_k $, 
we have $\hat{g}^{\infty}_{i}=((\Psi_j)^{-1}\circ \Psi_i)^{\ast} \hat{g}^{\infty}_{j}$ for the limit. Therefore,  there is a globally defined Riemannian metric $\hat{g}_{\infty}$ on $L$ such that its local expression in coordinates $\Psi_i:\Omega_i\rightarrow L$ is $\hat{g}_{\infty i}$.  Similarly,   the limits   $\bar{g}_{\infty i}$ and $X_{\infty i}$ on $\Omega_i$ also satisfy  the right change of variable formulas for different $i$,  and they also form two  globally defined tensor fields $\bar{g}_{\infty}$ and $X_{\infty}$ on $L$.

On each $\Omega_i$, the distance function of $\hat{g}_{\infty i}$ is the limit of the distance function of  $(\Psi^k_i)^{\ast} \hat{g}_k$. This implies that the distance function of $\hat{g}_{\infty}$ is the same as $d_{L}$ on each $\Psi_i(\Omega_i)$. Hence $d_{\hat{g}_{\infty}}=d_{L}$ holds globally. In particular, $\hat{g}_{\infty}$ is complete. 

By Theorems \ref{t4.1} and \ref{t4.2},    $\hat{g}_{\infty}$, $\bar{g}_{\infty}$ are smooth, flat and $X_{\infty}$ is smooth and parallel with respect to both $\hat{g}_{\infty}$ and  $\bar{g}_{\infty}$.

   On the other hand, since $\{\Psi^k_i(|z|\leq 1/3)\}_{i=1}^{n_l}$ is a covering of $B(\bar{w}_k, l) \subset {Q}_k$, by (\ref{5.29}), we have 
   $
   c^{-1} l^{4} \leq vol (B(\bar{w}_k, l)) \leq c n_l.
   $		
Because 	$\{\Psi_i(|z|\leq 1/40)\}_{i=1}^{n_l} \subset B_{L}(\bar{w}, l+1)$ are disjoint, we obtain 
$$
vol (B_{L}(\bar{w}, l+1)) \geq c^{-1} l^{4},  \ \ \text{for  all} \ \ l\geq 1.
$$
	That is to say,   $(L,\hat{g}_{\infty})$ is flat and has   maximal volume growth.  This implies   $(L,\hat{g}_{\infty})$  is isometric to $\mathbb{R}^4$.  Hence $(L, \bar{g}_{\infty}, X_{\infty})$	 is the standard  Minkowski spacetime $\mathbb{R}^{3,1}$	with  $X_{\infty}={\partial}/ {\partial t}$.

			Now we construct diffeomorphisms between large domains of ${Q}_k$ and $L$, which are close to $\Psi_i(\Psi^{k}_i)^{-1}\triangleq f^{k}_{i}$ on each $\Psi^{k}_i(|z|<1/3)$ in certain topologies. One can make use of the fact that $(L,\hat{g}_{\infty})$ is $\mathbb{R}^4$. Fix a partition of unity $\{ u^{k}_i\}_{i=1}^{n_k}$ subordinate to the covering $\{\Psi^k_i(|z|< 1/3)\triangleq U^{k}_i\}_{i=1}^{n_k}$ of $B(\bar{w}_k, k) \subset \bar{M}_k$ (i.e. $supp (u^k_i) \subset \subset U^{k}_i$) satisfying  
			$$
		1=\sum_{i=0}^{n_k} u^k_i(x)$$	
		for all $x\in B(\bar{w}_k, k)$. 	 Clearly, one can choose $u^{k}_i$ satisfying  
		$$
		||u^{k}_i\circ \Psi^k_i||_{W^{2,p_1}(|z|\leq 1/3)} \leq Cont..
		$$
		
		Define \begin{equation} \label{6.8}
		\begin{split}
		F_{k}: B(\bar{w}_k, k) &\rightarrow L\simeq \mathbb{R}^4 \\
		x & \mapsto \sum_{i=0}^{n_k} u^{k}_i(x) f^{k}_i(x),
		\end{split}
		\end{equation}
	where the summation  is  taken only over those indices $i$ with $U^{k}_i \owns x$.  For any fixed $i_0 \in \{0,1,\cdots n_k\}$, since $(\Psi^k_j)^{-1}\circ \Psi^k_{i_0} \overset{C^{1,\alpha}}{\rightarrow}  (\Psi_j)^{-1}\circ \Psi_{i_0}$ if $\Psi_j(\Omega_j) \cap \Psi_{i_0}(\Omega_{i_0})\neq \phi$, we have  
	\begin{equation}	\label{6.9}
	||(f^k_j-f^k_{i_0})\circ \Psi^{k}_{i_0}||_{C^{1,\alpha}}=||\Psi_j \circ  (\Psi^k_j)^{-1}\circ \Psi^{k}_{i_0}-\Psi_{i_0}||_{C^{1,\alpha}} \rightarrow 0,
	\end{equation}
which implies 	\begin{equation} \label{6.10}
||(F_k-f_{i_0}^{k}) \circ \Psi^{k}_{i_0}  ||_{C^{1,\alpha}(|z|\leq 1/3)}=||\sum^{n_k}_{i=0}u^{k}_i\circ \Psi^{k}_{i_0} \cdot (f^{k}_i-f^k_{i_0})\circ \Psi^{k}_{i_0}||_{C^{1,\alpha}(|z|\leq 1/3)}\rightarrow 0. 
\end{equation}

In particular, for any $l>0$, there is $k_0>0$ such that $F_{k}: B(\bar{w}_k,l) \rightarrow F_{k}(B(\bar{w}_k,l)) \subset L$ is a diffeomorphism for all $k\geq k_0$.  (\ref{6.10})  yields   that 
$$
|| (F^{-1}_k)^{\ast} \hat{g}_k-((f^{k}_{i_0})^{-1})^{\ast} \hat{g}_k||_{C^{\alpha}} \rightarrow 0 
$$
together with 
$$
|| ((f^{k}_{i_0})^{-1})^{\ast} \hat{g}_k-\hat{g}_{\infty}||_{C^{\alpha}} \rightarrow 0
$$
we have 
\begin{equation}\label{6.11}
|| (F^{-1}_k)^{\ast} \hat{g}_k-\hat{g}_{\infty}||_{C^{\alpha}} \rightarrow 0
\end{equation}
over all compact sets of $L$. 

Similarly,  by using  $||(F_k-f_{i_0}^{k}) \circ \Psi^{k}_{i_0}  ||_{C^{1,\alpha}(|z|\leq 1/3)}\rightarrow 0$, we have 
\begin{equation} \label{6.12}
||(F_k^{-1})^{\ast} \bar{g}_k-((f^{k}_i)^{-1})^{\ast} \bar{g}_k||_{C^{\alpha}} \rightarrow 0
\end{equation}
\begin{equation} \label{6.13}
|| (F_k^{-1})^{\ast} X_k-((f^{k}_i)^{-1})^{\ast} X_k||_{C^{\alpha}} \rightarrow 0
\end{equation}
hence 
\begin{equation} \label{6.14}
(F_k^{-1})^{\ast} \bar{g}_k\rightarrow \bar{g}_{\infty}, \ \ \ (F_k^{-1})^{\ast} X_k\rightarrow X_{\infty}\end{equation}
in $C^{\alpha}$ topology over all compact sets of $L$. The convergences in  (\ref{6.11}) and (\ref{6.14}) also take place  in weak $W^{1,p_1}$ topology. 
			
Let $y^{\alpha}$, $\alpha=0,1,2,3$, be standard coordinates on $L\simeq \mathbb{R}^{3,1}$	 such that 
\begin{equation} \label{6.15}
	\begin{split}	
\hat{g}_{\infty}&=(dy^0)^2+(dy^1)^2+(dy^2)^2+(dy^3)^2\\
\bar{g}_{\infty}&=-(dy^0)^2+(dy^1)^2+(dy^2)^2+(dy^3)^2\\
X_{\infty}&=\frac{\partial}{\partial y^0}.
\end{split}
					 \end{equation}	 
		 
Composing with the diffeomorphism $F_k$, $\{y^{\alpha}\}$ will be an approximating harmonic coordinate system for the metric $\hat{g}_k$. In the following, we shall perturb $\{y^{\alpha}\}$ a little so that it becomes $\hat{g}_k$-harmonic.  Fix a large $l\geq  100$, let $D_{l}=\{ y\in \mathbb{R}^4\simeq L :  |y|\leq l\}$, consider the Dirichlet problem on $D^{k}_l\triangleq F_k^{-1}(D_l)\subset {Q}_{k}$:
\begin{equation}	\label{6.16}
\begin{split}	
 \triangle_{\hat{g}_k} u^{\xi}  & =-\triangle_{\hat{g}_k} (y^{\xi} \circ F_k) \ \ \text{in}\   D^{k}_l\\
 u^{\xi} & =0 \ \ \text{on}\ \partial D^{k}_l\end{split}
\end{equation} 		
	where $\triangle_{\hat{g}_k}$ is the Laplace-Beltrami operator relative to $\hat{g}_k$. In coordinates $\{y^{\alpha}\}$, denote $h^{k}_{\alpha\beta}=((F^{-1}_k)^{\ast}\hat{g}_k)(\frac{\partial}{\partial y^{\alpha}}, \frac{\partial}{\partial y^{\beta}})$,   multiply  the  factor $\sqrt{det(h^k)}$ on both sides of (\ref{6.16}),  the equation  (\ref{6.16})  becomes of  divergence form:
	\begin{equation} \label{6.17}
	\begin{split}
	\partial_{\alpha} [(h^{k})^{\alpha\beta}\sqrt{det(h^k)} \partial_{\beta} u^{\xi}]&=-\partial_{\alpha} [(h^{k})^{\alpha\xi}\sqrt{det(h^k)} ]\\
	&=-\partial_{\alpha} [(h^{k})^{\alpha\xi}\sqrt{det(h^k)}-\delta^{\alpha\xi} ].		\end{split}
	\end{equation}

According to the  H$\ddot{o}$lder estimate for elliptic equations of divergence form (\ref{6.17}) (Theorem 8.24 and Theorem 8.27 in \cite{GT}): $\exists \alpha^{\prime}>0$, s.t.  
 \begin{equation} \label{6.18}
 ||u^{\xi}||_{C^{\alpha^{\prime}}(D_l)} \preceq \int_{D_l}| (h^{k})^{\alpha\xi}\sqrt{det(h^k)}-\delta^{\alpha\xi}|^q dy)^{\frac{1}{q}} \rightarrow 0
 \end{equation}
 	as $k\rightarrow \infty$, for fixed $q>n=4$.

	To prove that  $\{u^{\xi}+y^{\xi}\}$ is a coordinate system, we need $C^1$-estimate for $u^{\xi}$. For this purpose,  	 we  write the equation (\ref{6.16}) in harmonic coordinate system  	$\Psi^k_i:\{|z|<1/3\}\rightarrow {Q}_k$, where the Laplace operator $\triangle_{\hat{g}_k}$ has the form:
	\begin{equation} \label{6.19}
	\triangle_{\hat{g}_k} u^{\xi}=\hat{g}^{\alpha\beta}_{k} \partial_{\alpha}\partial_{\beta} u^{\xi}	\end{equation} 
	and $(\hat{g}_{k})_{\alpha\beta}=(\hat{g}_k)(\frac{\partial}{\partial z^{\alpha}}, \frac{\partial}{\partial z^{\beta}})$, $(\hat{g}_{k})^{\alpha\beta}=(\hat{g}^{-1}_{k})_{\alpha\beta}$. 
	From the formula 	\begin{equation} \label{6.20}
	\frac{\partial^2 y^{\xi}}{\partial z^{\alpha}\partial z^{\beta}}-\Gamma^{\delta}_{\alpha\beta}(\hat{g}_k)\frac{\partial y^{\xi}}{\partial z^{\delta}}+\Gamma^{\xi}_{\eta\sigma}({h^{k}})\frac{\partial y^{\eta}}{\partial z^{\alpha}}\frac{\partial y^{\sigma}}{\partial z^{\beta}}=0	,\end{equation}
and (\ref{6.16}), we have 
	\begin{equation} \label{6.21}
	||\triangle_{\hat{g}_k} u^{\xi}||_{L^{p_1}(|z|\leq 1/3)} \leq C.
	\end{equation}
	The  $L^{p}$ estimate further implies   
	\begin{equation} \label{6.22}
	||u^{\xi}||_{W^{2,p_1}(|z|\leq 1/3)}\leq C,
	\end{equation}
	since the coefficients of the equation (\ref{6.16}) (\ref{6.19})  satisfies $||\hat{g}^{\beta\gamma}_{k}-\delta_{\beta\gamma}||_{C^{\alpha}} \leq \epsilon_0$. By using 
	\begin{equation} \label{6.23}	\frac{\partial^2 u^{\xi}}{\partial y^{\alpha}\partial y^{\beta}}=\frac{\partial^2 u^{\xi}}{\partial z^{\delta}\partial z^{\eta}}\frac{\partial z^{\delta}}{\partial y^{\alpha}}\frac{\partial z^{\eta}}{\partial y^{\beta}}+\frac{\partial u^{\xi}}{\partial z^{\delta}}\frac{\partial^2 z^{\delta}}{\partial y^{\alpha}\partial y^{\beta}},	\end{equation}
	we know $||u^{\xi}||_{W^{2,p_1}_{y}(\Psi_i(|z|\leq 1/3))} \leq C$.  By Sobolev imbedding theorem, we have 
	\begin{equation} \label{6.24}
	||u^{\xi}||_{C^{1,\alpha_1}(|y|\leq l-1)}	\leq C.\end{equation}	
	
	Now we use an interpolation technique. For any $\delta>0$, $0<\alpha<\alpha_0$,  there exists $C=C(\delta, \alpha, \alpha_1, l)>0$ such that 
	\begin{equation} \label{6.25}
	||u^{\xi}||_{C^{1,\alpha}(|y|\leq l-1)}	\leq C ||u^{\xi}||_{C^{0}(|y|\leq l-1)}+\delta ||u^{\xi}||_{C^{1,\alpha_1}(|y|\leq l-1)},\end{equation}
	which  gives 	
	\begin{equation} \label{6.26}
	\limsup_{k\rightarrow \infty}||u^{\xi}||_{C^{1,\alpha}(|y|\leq l-1)}	\leq C \delta.
	\end{equation}	
		Since  $\delta$ is arbitrary, we obtain 
	\begin{equation} \label{6.27}	\lim_{k\rightarrow \infty}||u^{\xi}||_{C^{1,\alpha}(|y|\leq l-1)}=0.
	\end{equation}		
	
	Consequently, $\{x^{\xi}\triangleq y^{\xi}+u^{\xi}\}_{\xi=0}^{3}$	 is a $\hat{g}_k$-harmonic coordinate system on 
	$F^{-1}_k(|y|\leq l-1)$ for all large  $k$.

	Under the coordinate system $\{x^{\xi}\}$, we still have 
		\begin{equation} \label{6.28}
	\begin{split}
	& ||(\hat{g}_k)_{\beta\gamma}-\delta_{\beta\gamma}||_{C^{\alpha}(|x|\leq l-2)}\rightarrow 0\\
	& ||(\bar{g}_k)_{\beta\gamma}-\eta_{\beta\gamma}||_{C^{\alpha}(|x|\leq l-2)}\rightarrow 0\\
		& ||(X_k)^{\beta}-\delta_{\beta 0}||_{C^{\alpha}(|x|\leq l-2)}\rightarrow 0,\end{split}
	\end{equation}	
 where 
	$$
	(\hat{g}_{k})_{\beta\gamma}=(\hat{g}_{k})(\frac{\partial}{\partial x^{\beta}}, \frac{\partial}{\partial x^{\gamma}}), \ (\bar{g}_{k})_{\beta\gamma}=(\bar{g}_{k})(\frac{\partial}{\partial x^{\beta}}, \frac{\partial}{\partial x^{\gamma}}),\  (X_k)^{\beta}=dx^{\beta}(X_k).$$
	Note that we already have 
	\begin{equation} \label{6.29}
	\begin{split}
	 ||(\hat{g}_k)_{\beta\gamma}||_{W^{1,p_1}(|x|\leq \frac{l}{2})}+ ||(\bar{g}_k)_{\beta\gamma}||_{W^{1,p_1}(|x|\leq \frac{l}{2})}+ ||(X_k)^{\beta}||_{W^{1,p_1}(|x|\leq \frac{l}{2})} \leq C.\end{split}
	\end{equation}	
	We shall  prove in the final step	\begin{equation} \label{6.30}
	\begin{split}
	& ||(\hat{g}_k)_{\beta\gamma}-\delta_{\beta\gamma}||_{W^{1,p_1}(|x|\leq \frac{l}{2})}\rightarrow 0\\
	& ||(\bar{g}_k)_{\beta\gamma}-\eta_{\beta\gamma}||_{W^{1,p_1}(|x|\leq \frac{l}{2})}\rightarrow 0\\
		& ||(X_k)^{\beta}-\delta_{\beta 0}||_{W^{1,p_1}(|x|\leq \frac{l}{2})}\rightarrow 0,	\end{split}
	\end{equation}	
	which lead to  a contradiction with (\ref{6.1}). 
	
	To prove (\ref{6.30}),  we need to investigate the Einstein equations more carefully.  From  (\ref{3.9})-(\ref{3.14}) and (\ref{4.4})-(\ref{4.9}),  we have  
	
	 \begin{equation}  \label{6.31}	 \begin{split}
		\hat{g}^{\xi\eta}_{k} \partial_{\xi} \partial_{\eta} \log u_k & =f_1+\partial \tilde{f}_1	\end{split}
	\end{equation}	
	 \begin{equation}  \label{6.32}	 \begin{split}
	\hat{g}^{\xi\eta}_{k} \partial_{\xi} \partial_{\eta} (\hat{g}_k)& =f_2+\partial \tilde{f}_2	\\
	\end{split}
	\end{equation}
	and 
	\begin{equation} \label{6.33}
	\begin{split}
	d\Lambda_k& =f_3+\partial \tilde{f}_3\\
	\delta \Lambda_k&=f_4+\partial \tilde{f}_4
	\end{split}
	\end{equation}			
		where 
	\begin{equation} \label{6.34}
	\begin{split}
 |f_1|+|f_2|+|f_3|+|f_4|  \preceq & |Ric(\bar{g}_k)|+|\partial \hat{g}_k|^2+u^{-2}_k|\mathcal{L}_{X_k}\bar{g}_k|^2+|\partial \log u_k|^2\\& +u^{2}_k|\Lambda_k|^2\\
	 |\tilde{f}_1|+  |\tilde{f}_2| +|\tilde{f}_3|+|\tilde{f}_4| \preceq & u^{-1}_k|\mathcal{L}_{X_k}\bar{g}_k|.	\end{split}
	\end{equation}

	First we derive the $W^{1,2}$ estimate.   Let $\phi$ be a nonnegative smooth function  with compact support in $\{|x| \leq l-2\}$ such that $\phi\equiv 1$ on $\{|x|\leq l-3\}$ and $\phi\equiv 0$	on $\{l-2-\frac{1}{2}\leq |x|\leq l-2\}$.  From (\ref{6.31}),  we have  
	\begin{equation} \label{6.35}
	\begin{split}
	\int \phi^2 |\partial \log u_k|^2 &\leq \int\phi^2|\log u_k| (|u_k^{-1}\mathcal{L}_{{X}_{k}}\bar{g}_k|^2+u^{2}_k|\Lambda_k|^2+|Ric(\bar{g}_k)|) \\&  +\int |\partial \phi|^2|\log u_k|^2 +\int\phi^2(1+|\log u_k|^2) |u_k^{-1}\mathcal{L}_{{X}_{k}}\bar{g}_k|^2		\end{split}
	\end{equation}	
	where the formula (\ref{3.4}) for $D^{\hat{g}} X$ was used, which  yields 
	\begin{equation} \label{6.36}
	\int \phi^2 |\partial \log u_k|^2 \rightarrow 0.
	\end{equation}
	Multiplying  (\ref{3.13}) with $\phi^2$ and integrating  by parts, we have    
	\begin{equation}
	\begin{split} \label{6.37}
	&|-\int 2\phi \hat{g}^{\alpha\beta}_{k}\partial_{\alpha}\log u_k \partial_{\beta} \phi+\frac{1}{4}\int \phi^2 u_k^{2} |\Lambda_k|^2|\\
	&\leq \int |u_k^{-1}\mathcal{L}_{{X}_{k}}\bar{g}_k|(\phi |\partial \phi|+\phi^2|\partial \log u_k|+\phi^2 |u_k^{-1}\mathcal{L}_{{X}_{k}}\bar{g}_k|+\phi^2|Ric(\bar{g}_k)|)\\&\rightarrow 0,
	\end{split}
	\end{equation}
	which gives 
	\begin{equation} \label{6.38}
	\int \phi^2 u_k^{2} |\Lambda_k|^2 \rightarrow 0.
	\end{equation}	
	Similarly, from (\ref{6.32}), we have  
	\begin{equation} \label{6.39}
	\begin{split}
	 ||\phi\partial \hat{g}_k ||_{L^2}^2 \leq & ||(\hat{g}_k)_{\alpha\beta}-\delta_{\alpha\beta})||_{L^{\infty}}(1+||\phi u_k \Lambda_k||^2_{L^2}\\&+||u_k^{-1}\mathcal{L}_{{X}_{k}}\bar{g}_k||^2_{L^2}+||\phi \partial \log u_k||^2_{L^2})+||u_k^{-1}\mathcal{L}_{{X}_{k}}\bar{g}_k||^2_{L^2}\\
	 &\rightarrow 0.	 	\end{split}\end{equation}
	 
	 Consequently, 
	 \begin{equation} \label{6.40}
	 ||\partial \hat{g}_k ||_{L^2(|x|\leq l-3)}+ ||\partial \log u_k ||_{L^2(|x|\leq l-3)}+ || \Lambda_k ||_{L^2(|x|\leq l-3)} \rightarrow 0.	 \end{equation}

	By using the H$\ddot{o}$lder inequality,  for any $2<p<p_1$, we obtain   
	\begin{equation} \label{6.41}
	\begin{split}
	||\partial \hat{g}_k||_{L^{p}(|x|\leq l-3)} & \leq ||\partial \hat{g}_k||^{\frac{p_1}{p} \times \frac{p-2}{p_1-2}}_{L^{q}(|x|\leq l-3)}||\partial \hat{g}_k||^{\frac{2}{p}\times \frac{p_1-p}{p_1-2}}_{L^{2}(|x|\leq l-3)}	\rightarrow 0\\ ||\partial \log u_k||_{L^{p}(|x|\leq l-3)} & \leq ||\partial \log u_k||^{\frac{p_1}{p} \times \frac{p-2}{p_1-2}}_{L^{q}(|x|\leq l-3)}||\partial \log u_k||^{\frac{2}{p}\times \frac{p_1-p}{p_1-2}}_{L^{2}(|x|\leq l-3)}	\rightarrow 0 \\
	 || \Lambda_k||_{L^{p}(|x|\leq l-3)} & \leq ||\Lambda_k||^{\frac{p_1}{p} \times \frac{p-2}{p_1-2}}_{L^{p_1}(|x|\leq l-3)}||\Lambda_k||^{\frac{2}{p}\times \frac{p_1-p}{p_1-2}}_{L^{2}(|x|\leq l-3)}	\rightarrow 0.	\end{split}\end{equation}

	Before using  the $L^p$ estimate of second order elliptic equations with variable coefficients (see \cite{GT}),  we need to  treat the terms $\partial \tilde{f}_1$ and $\partial \tilde{f}_2$  in  (\ref{6.31}) and (\ref{6.32}). 	
	
Consider the Dirichlet problem 
	\begin{equation} \label{6.42}
	\begin{split}
	\hat{g}^{\xi\eta}_{k} \partial_{\xi} \partial_{\eta} F_1 & =\tilde{f}_1\ \ \ \text{in} \ \ \ \{|x|<l-3\}\\
	 F_1 & =0 \ \ \ \ \text{on} \ \ \ \{|x|=l-3\}.
	\end{split}
	\end{equation}	
 From (\ref{6.42}), we obtain    \begin{equation} \label{6.43}
	 ||\partial F_1||_{L^2}^2\preceq || \partial \hat{g}_k||^{\frac{1}{2}}_{L^4}||F_1||^{\frac{1}{2}}_{L^4}||\partial F_1||_{L^2} +||f_1||_{L^2} ||F_1||_{L^2}.	 \end{equation}
	 
	 Using Sobolev imbedding  $||F_1||_{L^4} \preceq  ||\partial F_1||_{L^2}$ and $|| \partial \hat{g}_k||_{L^4} \rightarrow 0$, we find  
	 \begin{equation} \label{6.44}
	 ||\partial F_1||_{L^2}^2\preceq ||\tilde{f}_1||_{L^2},  	 \end{equation}	 	
which particularly  implies that  the equation (\ref{6.42}) can be solved  in $W^{1,2}_0$.

	The interior $W^{2,p}$ estimate (Theorem 9.11 in \cite{GT}) on  (\ref{6.42}) gives 
	 \begin{equation} \begin{split}\label{6.45}
		|| {F}_1||_{W^{2,p} \{|x|<l-\frac{7}{2}\}}&\preceq ||\tilde{f}_1||_{L^{p}\{|x|<l-3\}}+|| {F}_1||_{W^{1,2} \{|x|<l-3\}} \\
		& \preceq ||\tilde{f}_1||_{L^{p}\{|x|<l-3\}} \end{split}	\end{equation}	
	for $p>2$. 	
		Note that  (\ref{6.31}) can be rewritten as 
			\begin{equation} \label{6.46}
		\hat{g}^{\xi\eta}_{k} \partial_{\xi} \partial_{\eta} (\log u_k-\partial F_1) =f_1+\hat{g}^{-2}_k\ast \partial \hat{g}_k \ast \partial^2F_1.		
			\end{equation}							
	By Sobolev imbedding theorem and interior $W^{2,p}$ estimate on  (\ref{6.46}), 
	\begin{equation} \begin{split} \label{6.47}
	& ||\partial (\log{u}_k-\partial F_1)||_{L^{p} (|x|\leq l-4)}\\ & \preceq  ||\log{u}_k-\partial F_1||_{W^{2, \frac{np}{n+p}} (|x|\leq l-4)} \\	
	& \preceq  ||f_1+\hat{g}^{-2}_k\ast \partial \hat{g}_k \ast \partial^2F_1||_{L^{\frac{np}{n+p}}(|x|\leq l-\frac{7}{2})}+||\log{u}_k-\partial F_1||_{L^{p} (|x|\leq l-\frac{7}{2})}	\end{split}
	\end{equation}	
	where $n=4<p_1$. 
	
The  H$\ddot{o}$lder inequality gives  $||\partial \hat{g}_k \ast \partial^2F_1||_{L^{\frac{np}{n+p}}} \preceq ||\partial \hat{g}_k||_{L^{n}}||\partial^2F_1||_{L^{p}}$.  For $p>2$,  it follows from  (\ref{6.45}) and (\ref{6.47}) that
		\begin{equation} \begin{split} \label{6.48}
	||\partial \log{u}_k||_{L^{p} (|x|\leq l-4)}  \preceq  & ||f_1||_{L^{\frac{np}{n+p}}(|x|\leq l-3)}+||\tilde{f}_1||_{L^{p}(|x|\leq l-3)}\\& +||\log{u}_k||_{L^{p} (|x|\leq l-3)}.	\end{split}
	\end{equation}
	Applying the same argument to (\ref{6.32}) and combining (\ref{6.48}),  we obtain 
	    \begin{equation} \begin{split} \label{6.49}
	& ||\partial \log{u}_k||_{L^{p_1} (|x|\leq l-4)}+ ||\partial \hat{g}_k||_{L^{p_1} (|x|\leq l-4)} \\ & \preceq ||\log{u}_k||_{L^{p_1} (|x|\leq l-3)} +||\hat{g}_k-\delta||_{L^{p_1} (|x|\leq l-3)}+ || u^{-1}_k \mathcal{L}_{X_k}\bar{g}_k ||_{L^{p_1}(|x|\leq l-3)}\\& \ \ \ +|| Ric(\bar{g}_k) ||_{L^{\frac{np_1}{n+p_1}}(|x|\leq l-3)}+||\partial \hat{g}_k||^2_{L^{\frac{2np_1}{n+p_1}} (|x|\leq l-3)}\\& \ \ \ +||\partial \log{u}_k||^2_{L^{\frac{2np_1}{n+p_1}} (|x|\leq l-3)}+|| u_k \Lambda_k||^2_{L^{\frac{2np_1}{n+p_1}} (|x|\leq l-3)} .\end{split}
	\end{equation}
	 Since  $\frac{2np_1}{n+p_1} <p_1$,  it follows from (\ref{6.30}),  (\ref{6.41}) and (\ref{6.49}) that  
		  \begin{equation} \begin{split} \label{6.50}
	||\partial \log{u}_k||_{L^{p_1} (|x|\leq l-4)}+ ||\partial \hat{g}_k||_{L^{p_1} (|x|\leq l-4)}\rightarrow 0.\end{split}
	\end{equation}	
	Now we are left with the estimate of $\Lambda_k$. In order to employ   the standard $L^p$ estimate of second order elliptic equations,  one can consider  the Dirichlet problem:
	\begin{equation} \label{6.51}
	\begin{split}
	 (d+\delta)^2 \eta& =\Lambda_k \ \ \ \text{in} \ \{|x|<l-4\}\\
	\eta & =0 \ \ \ \ \text{on} \ \ \{|x|=l-4\}.
	\end{split}
	\end{equation}
	
	 Recall  the well-known Bochner formula:
	 \begin{equation} \label{6.52}
	 \begin{split}
	 (d+\delta)^2\eta= & -\hat{g}_k^{\alpha\beta} \partial_{\alpha}\partial_{\beta}\eta+\hat{g}^{-1}_k \ast\partial \hat{g}_k \ast \partial {\eta}\\& + \partial (\hat{g}^{-1}_k \ast\partial \hat{g}_k \ast {\eta})+\hat{g}^{-2}_k\ast(\partial \hat{g}_k)^2\ast \eta
	 	 \end{split}
	 \end{equation}
	 where  the curvature term $Rm \ast \eta$ has been   rewritten as $\partial (\hat{g}^{-1}_k \ast \partial \hat{g}_k \ast {\eta})+\hat{g}^{-2}_k\ast (\partial \hat{g}_k)^2\ast \eta$. 	 
		 From (\ref{6.51}) and  (\ref{6.52}),  we have 
	\begin{equation} \label{6.53}
	 \begin{split}
	|| \partial  \eta||^2_{L^2} 
	& \preceq  ||\Lambda_k||_{L^2}||\eta||_{L^2}+ ||\partial \hat{g}_k||^2_{L^{4}} ||\eta||^{2}_{L^4}. 	 \end{split}
	 \end{equation}	
	As in (\ref{6.43}),  due to the  Sobolev inequality,  we find that  (\ref{6.51}) is solvable  and the solution $\eta$ satisfies 
	  \begin{equation} \label{6.54}
	  ||\eta||_{W^{1,2}(|x|\leq l-4)} \preceq ||\Lambda_k||_{L^2}\rightarrow 0.
	  \end{equation}
	 
	Furthermore,  taking  exterior derivative on  (\ref{6.51}),  we find  
	 \begin{equation} \label{6.55}
	(d+\delta)^2 \xi=d\Lambda_k
		\end{equation}	 
		where $\xi=d \eta$.  Analogous to (\ref{6.52}),  we rewrite the  equation (\ref{6.55})  as follows:
		\begin{equation} \label{6.56}
		\begin{split}
		-\hat{g}_k^{\alpha\beta} \partial_{\alpha}\partial_{\beta} \xi &=  \hat{g}_k^{-1}\ast \partial \hat{g}_k \ast \partial \xi+\partial (\hat{g}_k^{-1}\ast \partial \hat{g}_k \ast \xi)+ \hat{g}_k^{-2}\ast (\partial \hat{g}_k)^2 \ast \xi\\ &\ \  +f_3+\partial \tilde{f}_3\\
		&=f_5+\partial \tilde{f}_5
		\end{split}
		\end{equation}
		where $$f_5=f_3+ \hat{g}_k^{-1}\ast \partial \hat{g}_k \ast \partial \xi+ \hat{g}_k^{-2}\ast (\partial \hat{g}_k)^2 \ast \xi, \ \ \ \tilde{f}_5= \tilde{f}_3+\hat{g}_k^{-1}\ast \partial \hat{g}_k \ast \xi.$$  
		 From (\ref{6.56}) and direct computations,  
			\begin{equation}	 \label{6.57}
		 ||\partial(\phi \xi)||^2_{L^2} \preceq  || \xi||_{L^2}^2 +||\phi \xi||^2_{L^{\frac{2p_1}{p_1-2}}}  ||\partial \hat{g}_k ||^2_{L^{p_1}}+||f_3 ||^2_{L^2}+||\tilde{f}_3 ||^2_{L^2},
	\end{equation}	 
	where  $\phi$ is  a cutoff function on $\{|x|\geq l-9/4\}$. This  yields  
	 \begin{equation} \label{6.58}
		 ||\xi||_{W^{1,2}(|x|\leq l-5)} \rightarrow 0.
		 \end{equation}

		On the other hand,  applying  the same  estimate  (\ref{6.48}) to equation (\ref{6.56}),  it follows that for $q\leq p_1$, $5\leq i \leq l/3$  
		  		   \begin{equation} \begin{split}\label{6.59}
		 ||\partial \xi||_{L^{q}(|x|\leq l-i-1)} & \preceq  ||\tilde{f}_3||_{L^{q}(|x|\leq l-i)}+  ||f_3||_{L^{\frac{4q}{q+4}}(|x|\leq l-i)}+|| \xi||_{L^{q}(|x|\leq l-i)}\\	
		 &\ \ \ +||\partial \hat{g}_k||_{L^{p_1}(|x|\leq l-i)} ||\partial \xi||_{L^{\frac{1}{\frac{1}{q}-\frac{1}{p_1}+\frac{1}{4}}}(|x|\leq l-i)}\\
		 & \ \ \ + ||\partial \hat{g}_k||_{L^{p_1}(|x|\leq l-i)} || \xi||_{L^{\frac{qp_1}{p_1-q}}(|x|\leq l-i)}.		 \end{split}\end{equation}		
	By choosing $\frac{1}{q_1}=\frac{1}{2}-( \frac{1}{4}-\frac{1}{p_1})$ in (\ref{6.59}),  using (\ref{6.58}) and Sobolev imbedding, we have 
	\begin{equation} \label{6.60}	 ||\xi||_{W^{1, q_1}(|x|\leq l-6)} \rightarrow 0.\end{equation}
	Inductively,  choosing $\frac{1}{q_{i+1}}=\frac{1}{q_i}-( \frac{1}{4}-\frac{1}{p_1})$ if $q_{i+1} \leq p_1$, otherwise $q_{i+1}=p_1$ in (\ref{6.59}),  and repeating the process at most $[\frac{4p_1}{p_1-4}]$ times, 	we obtain 
	\begin{equation} \label{6.61}
		 ||d\eta||_{W^{1, p_1}(|x|\leq l-6-\frac{4p_1}{p_1-4})} \rightarrow 0.\end{equation}

	 The same argument applied  to the equation  
	 	$$
	(d+\delta)^2 \delta \eta=\delta \Lambda_k, 
	$$ 
      gives 
		  \begin{equation} \label{6.62}
	||\delta \eta||_{W^{1,p_1}(|x|\leq l-6-\frac{4p_1}{p_1-4})} \rightarrow 0.
	\end{equation}	 	
	
	Combining   (\ref{6.61}) (\ref{6.62}) (\ref{6.51}) and  (\ref{6.50}),  
	\begin{equation} \label{6.63}
	||\Lambda_k||_{L^{p_1} (|x|\leq l-6-\frac{4p_1}{p_1-4})}\rightarrow 0. 
		\end{equation}

Finally, combining (\ref{6.50})  (\ref{6.63}) and (\ref{6.30}), we have 
$$ r^{h}(W^{1,p_1}, \epsilon_1)(\bar{w}_k) \geq \frac{1}{2} l 
$$	
 for all $l>0$ when $k$ (possibly depending on $l$) is large enough.  
 This  leads to  a contradiction with (\ref{6.1}). The proof is completed. 
 \end{pf}
 
 \begin{pf} of Theorem \ref{t1.1}.
 
 The proof is  via an argument by contradiction proposed in Section 2.  Suppose the result is not true,   there is a  sequence of pointed Lorentzian manifolds with timelike vector fields  $({B}_{\hat{g}_j}(\bar{x}_j, \frac{1}{2}m_j), \bar{g}_j, \hat{g}_j, X_j, \bar{x}_j)$ satisfying (\ref{2.9}) and (\ref{2.10}). By Proposition \ref{p5.2},  for large $j$,   one can construct a non-degenerate map around $\bar{x}_j$ from a large Euclidean ball to the manifold ${B}_{\hat{g}_j}(\bar{x}_j, \frac{1}{2}m_j)$ such that the pull back of  $\hat{g}_j$ is comparable to the Euclidean metric (see (\ref{5.29})).  The contradiction is  derived ultimately from Proposition  \ref{p6.1}.
 \end{pf}
 \begin{pf} of Corollary \ref{c1.2}. 
 
 The estimate of $X$ has been incorporated in (\ref{2.4}) and proved finally  in Proposition \ref{p6.1}. 
 \end{pf}
 \begin{pf} of Corollary \ref{c1.3}. 
 
 Higher regularities  can be obtained by standard elliptic estimate on equations (\ref{3.9})-(\ref{3.14}) in $\hat{g}$-harmonic coordinates $\{z^{\alpha}\}$ established in Theorem \ref{t1.1}.   
 \end{pf} 
 
 \section{A generalization of Theorem \ref{t1.1}}
 First of all, we give a proof of Proposition \ref{p1.4}. The argument is motivated by \cite{CGT}. 
  \begin{pf} of Proposition  \ref{p1.4}.  
  
    Fix $\epsilon=10^{-4}$, and apply  Theorem \ref{t1.1} on  $\hat{B}(\bar{x}, \frac{a}{2})$ for any $\bar{x}\in\hat{B}({x}_0, \frac{a}{2})$, let $\Psi^{\bar{x}}: \{|z|< \frac{1}{2}c^{-1}_1{a}\}\rightarrow \bar{M}$ be the immersion map with  $\Psi^{\bar{x}}(\{ |z|< \frac{1}{2}c^{-1}_1{a}\}) \supset \hat{B}(\bar{x}, \frac{1}{2}c_2^{-1} a)$.    For any point $x\in \hat{B}(\bar{x}, \frac{1}{4}c_2^{-1} a)$,  we claim that the number of the pre-images  $(\Psi^{\bar{x}})^{-1}(x) \subset \{ |z|< \frac{1}{4}c^{-1}_1{a}\}$ satisfies 
\begin{equation} \label{7.1}
 {c}^{-1} \frac{a^4}{vol_{\hat{g}}(\hat{B}(\bar{x}, \frac{1}{4}c^{-1}_2{a}))}  \leq \# (\Psi^{\bar{x}})^{-1}(x)\leq {c} \frac{a^4}{vol_{\hat{g}}(\hat{B}(\bar{x}, \frac{1}{4}c^{-1}_2{a}))}.
\end{equation}
To prove (\ref{7.1}), let $x_1,x_2,\cdots, x_m$ be all the pre-images of ${x}$ in $\{|z|<\frac{1}{4}c^{-1}_1{a}\}$. For any point $\tilde{x}\in \hat{B}(\bar{x}, \frac{1}{4}c^{-1}_2{a}))$, let  $\gamma:[0,b]\rightarrow \hat{B}(\bar{x}, \frac{1}{2}c^{-1}_2{a})$ be a $\hat{g}$-geodesic connecting ${x}$ and $\tilde{x}$ with length $b \leq \frac{1}{2}c^{-1}_2{a}$. Let $\gamma_i$ be the lifting of $\gamma$  with $\gamma_i(0)=x_i$. Then $\{\gamma_i(b)\}$ will be $m$-distinct pre-images of $\tilde{x}$ in $\{|z|<\frac{1}{2}c^{-1}_1{a}\}$. This implies 
\begin{equation} \label{7.2}
m \cdot vol_{\hat{g}}(\hat{B}(\bar{x}, \frac{1}{4}c^{-1}_2{a}) \leq vol_{\hat{g}}(\{|z|<\frac{1}{2}c^{-1}_1{a}\}).
\end{equation}

On the other hand,  let   $\tilde{x}_1,\tilde{x}_2,\cdots, \tilde{x}_{\tilde{m}}$ be all pre-images of $\tilde{x}$ in $\{ |z|< \frac{1}{8}c^{-1}_1{a}\}$.  The above argument gives  $\tilde{m}\leq m$. Then 
\begin{equation}  \label{7.3} vol_{\hat{g}}(\{|z|<\frac{1}{8}c^{-1}_1{a}\})  \leq  m\cdot vol_{\hat{g}}(\hat{B}(\bar{x}, \frac{1}{4}c^{-1}_2{a}).
\end{equation}

(\ref{7.1}) follows from (\ref{7.2}) and (\ref{7.3}).

By Theorem \ref{t5.1} and (\ref{2.5}), there is a Riemannian metric $\check{g}$ on $B_{\hat{g}}(x_0, \frac{3}{4}a)$ with $\frac{1}{2}\hat{g}\leq \check{g}\leq 2 \hat{g}$ and $|Rm(\check{g})| \leq ca^{-2}$.  By using the  volume comparison theorem (for $\check{g}$) in Riemannian geometry, we have 
\begin{equation} \label{7.4}
vol_{\hat{g}}(\hat{B}(\bar{x}, \frac{1}{4}c^{-1}_2{a})) \geq c^{-1}vol_{\hat{g}}(\hat{B}(x_0, \frac{1}{2} a)).
\end{equation}
Combining 
(\ref{7.1}), (\ref{7.4}) and (\ref{1.11}), it  yields  
\begin{equation} \label{7.5}
 {c}^{-1} v^{-1}_0 \leq \# (\Psi^{\bar{x}})^{-1}(x)\leq {c} v^{-1}_0.
 \end{equation}
We claim  that for any two distinct pre-images $x_1,x_2\in (\Psi^{\bar{x}})^{-1}(x)$, 
\begin{equation} \label{7.6}
d_{\hat{g}}(x_1,x_2) \geq c^{-1} v_0 a
\end{equation}
which  will imply  that  the map $\Psi^{\bar{x}}$ is injective on $\{|z|\leq c^{-1} v_0 a\}$. 

Note that the metric $\check{g}$ has bounded curvature ($\leq ca^{-2}$) and volume bounded from below, i.e.  $vol_{\check{g}}(\hat{B}(x_0, \frac{3}{4}a))\geq c^{-1} v_0 a^4$. By  Theorem 4.3 in \cite{CGT},  the injectivity radius of $\check{g}$ at any point $x\in \hat{B}(x_0, \frac{1}{2}a)$ is bounded from below by $c^{-1} v_0 a$. If (\ref{7.6}) were not true for $\hat{g}$, it is also not true for $\check{g}$,  draw a shortest $\check{g}$-geodesic $\tilde{\gamma}$ connecting $x_1$ and $x_2$, the projection $\gamma=\Psi^{\bar{x}}(\tilde{\gamma})$ will be  a  closed and short $\check{g}$-geodesic at $x$, leading to  a contradiction with the injectivity radius estimate. 
\end{pf}

 We discuss a  generalization of Theorem \ref{t1.1}.  The proof also uses  the spirit of covering maps as in  Proposition \ref{p1.4}.  The point is to observe that,   in the proof of Theorem \ref{t1.1},   instead of the $L^{\infty}$ bounds,  only some uniform  normalized   $L^{p}$ bounds   of   $\mathcal{L}_{X}\bar{g}$ (for $p>4$) and  $Ric(\bar{g})$ (for $p>2$)  are needed.   
  \begin{thm} \label{t7.1}  For any  $\epsilon>0$,   $p_1>4$, $p_2>2$, there is a constant $\delta=\delta(p_1, p_2, \epsilon)>0$ with the following properties. Suppose a Lorentzian manifold $(\bar{M}, \bar{g})$ with a timelike vector field $X$ satisfies \begin{equation} \label{7.7}
 r (\frac{1}{vol(B_{\hat{g}}(x,r))}\int_{B_{\hat{g}}(x,r)}|u^{-1}\mathcal{L}_X\bar{g}|_{\hat{g}}^{p_1}dv_{\hat{g}})^{\frac{1}{p_1}} \leq \delta,
\end{equation}
and 
\begin{equation} \label{7.8}
 r^2 (\frac{1}{vol(B_{\hat{g}}(x,r))}\int_{B_{\hat{g}}(x,r)}|Ric(\bar{g})|_{\hat{g}}^{p_2}dv_{\hat{g}})^{\frac{1}{p_2}} \leq \delta, 
\end{equation}
for all $x\in B_{\hat{g}}(x_0,a)$ and  $0<r<d(x, \partial B_{\hat{g}}(x_0,a))$. 
Then the conclusion of Theorem \ref{t1.1} holds for $p=\min\{ p_1, \frac{4p_2}{4-p_2}\}$. 
\end{thm}
\begin{pf}
 We mimic the proof of Theorem \ref{t1.1} from Section 2 to Section 6. Suppose (\ref{2.5}) is not true, there are small   $\epsilon_0>0$  and a sequence of pointed Lorentzian manifolds $(\bar{M}_j, \bar{g}_j, X_j, P_j)$ with timelike vector fields   satisfying    (\ref{7.7}) (\ref{7.8}) for  
  $\delta_j\rightarrow 0$, but  (\ref{2.7}) holds  for  $m_j\rightarrow \infty$.

  The point is that the conditions  (\ref{7.7}) and (\ref{7.8}) are scaling invariant.  This implies that  for the sequence of scaled metrics in Section 2, we still have   
     \begin{equation} \label{7.9}
 r (\frac{1}{vol(B_{\hat{g}_j}(x,r))}\int_{B_{\hat{g}_j}(x,r)}|u_j^{-1}\mathcal{L}_{X_j}\bar{g}_j|_{\hat{g}_j}^{p_1}dv_{\hat{g}_j})^{\frac{1}{p_1}} \leq \delta_j
\end{equation}
and 
\begin{equation} \label{7.10}
 r^2 (\frac{1}{vol(B_{\hat{g}_j}(x,r))}\int_{B_{\hat{g}_j}(x,r)}|Ric(\bar{g}_j)|_{\hat{g}_j}^{p_2}dv_{\hat{g}_j})^{\frac{1}{p_2}} \leq \delta_j
 \end{equation}  
 for all   $x\in B_{\hat{g}_j}(\bar{x}_j,\frac{1}{2}m_j)$ and  
 $ 0<r\leq d(x, \partial B_{\hat{g}_j}(\bar{x}_j,\frac{1}{2}m_j))$. Another  observation is that  (\ref{7.9}) and (\ref{7.10}) remain almost invariant  for the pulled back metric under a good local covering map (replacing $\delta_j$ with $C\delta_j$).
By virtue of (\ref{7.1}),  the conditions  (\ref{7.9}) and (\ref{7.10}) imply that  on the local immersion chart $\Psi:\{z\in \mathbb{R}^4: |z|< \frac{1}{3}\}\rightarrow \bar{M}_j$ in Section 4.1,  we have  
   \begin{equation} \label{7.11}
 \int_{\{|z|\leq \frac{1}{2}\}}|u_j^{-1}\mathcal{L}_{X_j}\bar{g}_j|_{\hat{g}_j}^{p_1}dv_{\hat{g}_j}+ \int_{\{|z|\leq \frac{1}{2}\}}|Ric(\bar{g}_j)|_{\hat{g}_j}^{p_2}dv_{\hat{g}_j} \leq C\delta_j \rightarrow 0.\end{equation}
 
 It can be checked that  under (\ref{7.11}) the argument in Section 4.1 can go through for $p=\min\{p_1, \frac{4p_2}{4-p_2}\}>4$. Moreover, under the local covering map $\Phi_j$ in Proposition \ref{p5.2}, the conditions  (\ref{7.9}) and (\ref{7.10}) hold for $(Q_j, \Phi^{\ast}_j \bar{g}_j, \Phi^{\ast}_j \hat{g}_j, \Phi^{\ast}_j X_j)$  with $C\delta_j$ replacing $\delta_j$, where $Q_j=\{w: d_{\hat{g}_j}(w,0)<\frac{a_j}{2}\}$. This ensures that the argument ((\ref{6.49}) and (\ref{6.59}))  in Section 6 can also go through for $p=\min\{p_1, \frac{4p_2}{4-p_2}\}$.  
\end{pf}

 \textbf{Acknowledgement } 
 The work was partially supported  by grants 2022YFA1005400,  NSFC12141106.

Bing-Long Chen \\
 Department of Mathematics, \\
Sun Yat-sen University,\\
Guangzhou, P.R.China, 510275\\
Email: mcscbl@mail.sysu.edu.cn\\
\begin{thebibliography}{99}

 \bibitem{A1} M. T.  Anderson,   {\sl On stationary vacuum solutions to the Einstein equations}, Ann. H. Poincar$\acute{e}$    \textbf{1} (2000), no. \textbf{5}, 977-994.
 \bibitem{A2} M. T.  Anderson,   {\sl On long-time evolution in general relativity and geometrization of 3-manifolds},  Comm. Math. Phys. \textbf{222} (2001), 533-567.

\bibitem{CGT} J. Cheeger, M. Gromov, M. Taylor, {\sl  Finite propagation speed, kernel estimates for functions of the
Laplace operator, and the geometry of complete Riemannian manifolds}. J. Diff. Geom. \textbf{17} (1982),  15-53.

\bibitem{Chen1} B. L. Chen,  {\sl On stationary solutions to the vacuum Einstein field equations},  Asian J.  Math. \textbf{23} (2019), no. 4, 609-630.

 \bibitem{CB} Y. Choquet-Bruhat, {\sl Espace-temps Einsteiniens generaux, chocs gravitationnels}, Ann. Inst. Henri Poincar$\acute{e}$, \textbf{8}, 327-38, 1968.
  
\bibitem{EP} A. Einstein, W. Pauli, {\sl On the nonexistence of regular stationary solutions to relativistic field equations}, Ann. of Math. \textbf{44} (1943), 131-137. 

\bibitem{GT} D. Gilbarg, N. Trudinger, {\sl Elliptic Partial Differential Equations of Second Order}, Classics in Mathematics, ISBN 7-5062-5922-2, Springer. 


\bibitem{H} A. Haefliger, {\sl Groupoids and Foliations},  in: Groupoids in Analysis, Geometry, and Physics,  Contemp. Math. \textbf{282}, Am. Math. Soc., 83-100 (2001).


\bibitem{HE} S. Hawking, G.F.R. Ellis,  {\sl The Large Scale Structure of Spacetime},  Cambridge University Press,  1973. 


\bibitem{JK} J. Jost,  H.  Karcher,  {\sl Geometrische methoden zur gewinnung von a-priori-schranken f$\ddot{u}$r harmonische Abbildungen}, Manu. Math. \textbf{40} (1982), no.1, 27-77. 

\bibitem{KR10} S. Klainerman, I. Rodnianski, {\sl On the breakdown criterion in general relativity},  J. Amer. Math. Soc. \textbf{23} (2010), no. 2, 345-382.


\bibitem{KRS} S. Klainerman, I. Rodnianski, J. Szeftel, {\sl On bounded $L^2$ curvature conjecture},  Invent. Math. (2015) \textbf{202}, 91-216. 




\bibitem{Lich55} A. Lichnerowicz,  {\sl Theories Relativistes de la Gravitation et de L' Electromagnetisme},  Masson  and  Cie.,  Paris,  1955. 

\bibitem{Lich67} A. Lichnerowicz,  {\sl Relativistic Hydrodynamics and Magnetohydrodynamics},  W. A. Benjamin, Inc. New Yok, 1967. 


\bibitem{Lott} J. Lott,  {\sl On the long time behavior of Type-III Ricci flow solutions}, Math. Ann. \textbf{339} (2007): 627-666. 





\bibitem{Pe} P. Petersen,   {\sl Riemannian Geometry},  GTM \textbf{171},  Springer.  

\bibitem{PWY} P. Petersen,  G.F. Wei,  R.G. Ye,  {\sl Controlled geometry via smoothing}, Comm. Math. Helv. \textbf{74} (1999), 345-363. 

\bibitem{W} Q. Wang, {\sl Improved breakdown criterion for Einstein vacuum equations in cmc gauge},  Comm. Pure Appl. Math., Vol. LXV, 0021-0076 (2012).   
  \end{thebibliography}
\end{document}